\documentclass[12pt]{article}

\usepackage[margin=1in]{geometry}
\usepackage{amsthm, amssymb, amsmath, amsfonts}
\usepackage{mathrsfs}

\newtheorem{theorem}{Theorem}[section]
\newtheorem{lemma}[theorem]{Lemma}
\newtheorem{proposition}[theorem]{Proposition}
\newtheorem{corollary}[theorem]{Corollary}
\theoremstyle{definition}
\newtheorem{definition}[theorem]{Definition}
\newtheorem{example}[theorem]{Example}
\newtheorem{remark}[theorem]{Remark}
\newtheorem{question}[theorem]{Question}

\newcommand{\mibb}[1]{\left[#1\right]}
\newcommand{\bK}{\mathbb{K}}
\newcommand{\bD}{\mathbb{D}}
\newcommand{\bN}{\mathbb{N}}
\newcommand{\bR}{\mathbb{R}}
\newcommand{\bC}{\mathbb{C}}
\newcommand{\cA}{\mathcal{A}}
\newcommand{\cL}{\mathscr{L}}
\newcommand{\cR}{\mathscr{R}}
\newcommand{\cS}{\mathcal{S}}
\newcommand{\cN}{\mathcal{N}}
\newcommand{\cB}{\mathcal{B}}
\newcommand{\DR}{\mathscr{D}\!\mathscr{R}}
\newcommand{\MR}{\mathscr{M}\!\mathscr{R}}
\newcommand{\EDR}{\mathscr{E}\!\mathscr{D}\!\mathscr{R}}
\newcommand{\EMR}{\mathscr{E}\!\mathscr{M}\!\mathscr{R}}
\newcommand{\Fh}{\widehat F}
\newcommand{\Fb}{\overline{\widehat F}}
\newcommand{\pb}{\overline{\varphi}}
\newcommand{\fb}{\overline{f_{12}}}
\newcommand{\fr}[1]{f_{[#1]}}
\newcommand{\Psr}[1]{\Psi_{[#1]}}
\newcommand{\ue}{\mathrm{e}}
\newcommand{\ff}[2]{#1^{\underline{#2}}}
\newcommand{\ev}{{\rm ev}}
\newcommand{\od}{{\rm od}}

\allowdisplaybreaks

\begin{document}
\sloppy

\title{Exponential Double and Multiple Riordan Arrays\\ as Fiber Products}

\author{Tian-Xiao He\\
{\small Department of Mathematics}\\
{\small Illinois Wesleyan University}\\
{\small Bloomington, IL 61702-2900, USA}\\
}
\date{}
\maketitle

\begin{abstract}
\noindent
A double Riordan array has nonzero entries only in positions $(n,k)$ with $n+k$ even, and its even and odd parts are Riordan arrays with the same multiplier function. We develop this observation systematically, with the emphasis on exponential double and exponential multiple Riordan arrays. Using generalized Riordan arrays with respect to a sequence $(c_n)$, we show that the double Riordan group, the exponential double Riordan group, the multiple Riordan group, and the exponential multiple Riordan group are fiber products of (generalized) Riordan groups over the Lagrange group. For exponential double and exponential multiple Riordan arrays this yields product and inverse formulas without square or $\ell$th roots, valid over any field of characteristic zero; an abelian normal subgroup with quotient the Lagrange group and a corresponding splitting; subgroups isomorphic to the Riordan group; an iterated split-extension description of the multiple groups; and a characterization of elements of finite order, including the fact that over ${\mathbb R}$ every element of finite order is an involution. We show that the $\ell$th power of the production matrix is the direct sum of the production matrices of the compressions, give an explicit formula for it in terms of partial Bell polynomials, and characterize these arrays by a shared $A$-sequence. Finally, we reverse the construction. Two Riordan arrays with a common multiplier can always be interleaved into a unique double Riordan array, and into an ordinary Riordan array exactly when their first functions satisfy one more condition. Examples include Pascal-, Catalan-, Morgan--Voyce-, and ballot-type arrays.


\end{abstract}

\section{Introduction}

Riordan arrays and the Riordan group were introduced by Shapiro, Getu, Woan, and Woodson \cite{SGWW}. Davenport, Shapiro, and Woodson \cite{DSW} introduced double Riordan arrays $(g;f_1,f_2)$, whose columns are generated by an even function $g$ and two odd multipliers $f_1$ and $f_2$ used alternately; their sequence characterizations and compressions were studied in \cite{He18}, and their $A$- and $Z$-sequences in \cite{BDFJT, DFSW}. Triple Riordan arrays were considered by Barry \cite{Barry}, and the multiple Riordan group with $\ell$ multipliers in \cite{He-MR}. Double almost-Riordan arrays and their group were studied in \cite{He-DAR}, and their sequence characterization, compressions, and total positivity in \cite{He-DARA}. Exponential double and multiple Riordan arrays were introduced in \cite{He-EDRA, He-EMRA}. In all of these, the group law involves the square root $\sqrt{f_1f_2}$, or the $\ell$th root $(f_1\cdots f_\ell)^{1/\ell}$; over $\bR$ this imposes sign conditions, and in general it hides the structure of the groups.

The purpose of this paper is to explain these groups, in particular the exponential double and exponential multiple Riordan groups, by a classical construction from group theory.

\begin{definition}\label{def:fiber}
Let $\varphi_i:G_i\to H$, $1\leq i\leq m$, be group homomorphisms. The ($m$-fold) {\it fiber product} of the $G_i$ over $H$ is the subgroup
\[
G_1\times_HG_2\times_H\cdots\times_HG_m=\{(x_1,\ldots,x_m)\in G_1\times\cdots\times G_m:\ \varphi_1(x_1)=\cdots=\varphi_m(x_m)\}
\]
of the direct product, with componentwise multiplication.
\end{definition}

It is a subgroup because $\varphi_i(x_ix_i')=\varphi_i(x_i)\varphi_i(x_i')$. If every $\varphi_i$ is surjective with kernel $K_i$, then
\begin{equation}\label{exact}
1\longrightarrow K_1\times\cdots\times K_m\longrightarrow G_1\times_H\cdots\times_HG_m\longrightarrow H\longrightarrow1
\end{equation}
is exact, where the last map is $(x_1,\ldots,x_m)\mapsto\varphi_1(x_1)$. For example, if $\varphi_1=\varphi_2={\rm sgn}:S_3\to\{\pm1\}$, the fiber product $S_3\times_{\{\pm1\}}S_3$ consists of the $18$ pairs of permutations with the same sign.

The starting observation is that a double Riordan array is a {\it checkerboard} matrix, i.e., its entries vanish in positions $(n,k)$ with $n+k$ odd. It is therefore determined by its even part $(d_{2m,2k})$ and its odd part $(d_{2m+1,2k+1})$, and these {\it compressions} are Riordan arrays with the same multiplier function. Labelling a Riordan array $(G,F)$ by its multiplier $F$, the double Riordan group is the fiber product $\cR\times_\cL\cR$ of two copies of the Riordan group over the Lagrange group $\cL$. For exponential arrays, the compressions are generalized Riordan arrays with respect to the sequences $((2n)!)$ and $((2n+1)!)$, or $((\ell n+r)!)$ in the multiple case, which all lie over the same Lagrange group.

We now describe the main results. Let $\cR_c$ denote the generalized Riordan group with respect to $c=(c_n)$ (Section \ref{sec:gen}).
\begin{itemize}
\item The double Riordan group satisfies $\DR\cong\cR\times_\cL\cR$ (Theorem \ref{thm:DR}), and the exponential double Riordan group satisfies $\EDR\cong\cR_{(2n)!}\times_\cL\cR_{(2n+1)!}$ (Theorem \ref{thm:main2}). This yields root-free product and inverse formulas (Corollary \ref{cor:rootfree2}), an abelian normal subgroup $\cN\cong(\bK[[x]]^\times)^2$ with $\EDR\cong\cN\rtimes\cL$ (Theorem \ref{thm:normal2}), a characterization of involutions (Theorem \ref{thm:involution}), the identity $P^2=P_{D^\ev}\oplus P_{D^\od}$ for the production matrix (Theorem \ref{thm:P2}) with an explicit formula (Theorem \ref{thm:P2explicit}), and a shared $A$-sequence characterization (Theorem \ref{thm:sharedA2}).
\item For $\ell\geq2$, the exponential multiple Riordan group $\EMR_\ell$ and the (ordinary) multiple Riordan group are $\ell$-fold fiber products (Theorem \ref{thm:mainL}). The associated subgroup of $\EMR_\ell$ is isomorphic to $\EMR_{\ell-1}$, so $\EMR_\ell$ is obtained from the Riordan group by $\ell-1$ successive split extensions by the Appell group (Theorem \ref{thm:iterated}). Elements of finite order are characterized, and over $\bR$ they are involutions (Theorem \ref{thm:finite}, Corollary \ref{cor:real}). The production matrix has a cyclic block structure, $P^\ell$ is the direct sum of the production matrices of the $\ell$ compressions, and $P^\ell$ is given explicitly by partial Bell polynomials (Theorems \ref{thm:blocks}, \ref{thm:Plexplicit}).
\item Conversely (Section \ref{sec:inter}), two Riordan arrays $X=(G_0,F_0)$ and $Y=(G_1,F_1)$ interleave into a double Riordan array if and only if $F_0=F_1$, and into an ordinary Riordan array if and only if, in addition, $F_0=x(G_1/G_0)^2$ (Theorems \ref{thm:interleave-double} and \ref{thm:interleave-riordan}).
\end{itemize}

The paper is organized as follows. Section \ref{sec:gen} recalls generalized and exponential Riordan arrays and the properties used later, and proves a lemma on periodic matrices and their compressions. Section \ref{sec:EDR} treats exponential double Riordan arrays in detail, after a brief account of the ordinary double Riordan arrays, which supplies the coordinates used throughout. Section \ref{sec:EMR} treats exponential multiple Riordan arrays together with their ordinary counterparts. Section \ref{sec:inter} is devoted to interleaving, and Section \ref{sec:questions} lists some questions.

Throughout, $\bK$ is a field of characteristic zero (for example $\bR$ or $\bC$). For ordinary arrays, everything holds over a commutative ring $\bD$ with ``nonzero'' replaced by ``invertible''. All matrices are infinite and lower triangular with rows and columns indexed by $\bN_0=\{0,1,2,\ldots\}$, so all products involve only finite sums. For a lower triangular matrix $X$ we write $\overline X=UX$ for $X$ with its first row deleted, where $U=(\delta_{i+1,j})_{i,j\geq0}$ is the upper shift matrix, and $P_X=X^{-1}\overline X$ for its production matrix.

\section{Generalized and exponential Riordan arrays}\label{sec:gen}

Let $c=(c_n)_{n\geq0}$ be a sequence of nonzero elements of $\bK$, and let $\Lambda_c={\rm diag}(c_0,c_1,c_2,\ldots)$.

\begin{definition}\label{def:gen}
For $G\in\bK[[x]]$ with $G(0)\not=0$ and $F\in\bK[[x]]$ with $F(0)=0\not=F'(0)$, the {\it generalized Riordan array with respect to $c$} is the matrix $(G,F)_c=(d_{n,k})_{n,k\geq0}$ with
\begin{equation}\label{gen.def}
d_{n,k}=\mibb{\frac{x^n}{c_n}}G(x)\frac{F(x)^k}{c_k}=\frac{c_n}{c_k}\,[x^n]\,G(x)F(x)^k,
\end{equation}
i.e., the $k$th column of $(G,F)_c$ has the $c$-generating function $GF^k/c_k$. For $c_n=1$ this is the ordinary Riordan array $(G,F)$, and for $c_n=n!$ it is the exponential Riordan array, which we denote by $[G,F]$.
\end{definition}

Equation \eqref{gen.def} says that
\begin{equation}\label{gen.conj}
(G,F)_c=\Lambda_c\,(G,F)\,\Lambda_c^{-1}.
\end{equation}

\begin{proposition}\label{prop:gen}
{\rm (i) (Fundamental theorem)} If $b(x)=\sum_kb_kx^k/c_k$ and $b$ is identified with the column vector $(b_0,b_1,\ldots)^T$, then $(G,F)_c\,b=G\cdot b(F)$, as $c$-generating functions.

{\rm (ii)} The set $\cR_c$ of generalized Riordan arrays with respect to $c$ is a group under matrix multiplication, with
\begin{equation}\label{gen.prod}
(G,F)_c(D,H)_c=\big(G\,D(F),\ H(F)\big)_c,\qquad (G,F)_c^{-1}=\Big(\frac{1}{G(\overline F)},\ \overline F\Big)_c,
\end{equation}
where $\overline F$ is the compositional inverse of $F$, and $(G,F)_c\mapsto\Lambda_c^{-1}(G,F)_c\Lambda_c=(G,F)$ is an isomorphism $\cR_c\cong\cR$ onto the Riordan group $\cR=\cR_{(1)}$.
\end{proposition}

\begin{proof}
(i) $\sum_kd_{n,k}b_k=[x^n/c_n]\sum_kGb_kF^k/c_k=[x^n/c_n]\,G\,b(F)$. (ii) follows from \eqref{gen.conj} and the corresponding facts for $\cR$ \cite{SGWW}.
\end{proof}

Let $\cL=\{F\in\bK[[x]]:\ F(0)=0\not=F'(0)\}$ with the operation $F\star H=H\circ F$; it is a group (the {\it Lagrange group}) and is identified with the Lagrange subgroup $\{(1,F)\}$ of $\cR$. Let $\bK[[x]]^\times=\{G\in\bK[[x]]:\ G(0)\not=0\}$ be the abelian group of invertible power series under multiplication.

\begin{proposition}\label{prop:pi}
For every $c$, the label map
\[
\pi_c:\cR_c\to\cL,\qquad \pi_c\big((G,F)_c\big)=F,
\]
is a surjective homomorphism. Its kernel is the Appell subgroup $\cA_c=\{(G,x)_c\}$, and $G\mapsto(G,x)_c$ is an isomorphism $\bK[[x]]^\times\cong\cA_c$. Moreover, the isomorphisms $\cR_c\cong\cR$ of Proposition \ref{prop:gen} are compatible with the label maps.
\end{proposition}

\begin{proof}
By \eqref{gen.prod}, the multiplier of a product is $H\circ F=F\star H$, so $\pi_c$ is a homomorphism; it is onto since $\pi_c((1,F)_c)=F$. Its kernel consists of the arrays with $F=x$. Since $(G,x)_c(D,x)_c=(GD,x)_c$ and $(G,x)_c$ determines $G$ (its first column is $(c_n/c_0)[x^n]G$), $G\mapsto(G,x)_c$ is an isomorphism. Conjugation by $\Lambda_c$ does not change $F$.
\end{proof}

Changing the scale $c$ rescales the entries but not the label $F$; this is why generalized Riordan groups with different scales can be combined in a fiber product over $\cL$. By Proposition \ref{prop:pi}, all fiber products $\cR_{c^{(0)}}\times_\cL\cdots\times_\cL\cR_{c^{(m-1)}}$ are isomorphic to $\cR\times_\cL\cdots\times_\cL\cR$.

For $c_n=n!$ we write $[g,f]=(g,f)_{n!}$, so that $d_{n,k}=\frac{n!}{k!}[t^n]g(t)f(t)^k$ and the $k$th column has the exponential generating function $gf^k/k!$. We use the following facts \cite{DFR, He15}.
\begin{itemize}
\item[(E1)] $[g,f]\,h(t)=g(t)\,h(f(t))$ for exponential generating functions $h$; in particular $[g,f]^{-1}B=B(\bar f)/g(\bar f)$, where $\bar f$ is the compositional inverse of $f$.
\item[(E2)] Deleting the first row of a column corresponds to differentiating its exponential generating function. Hence $\overline{[g,f]}$ has columns $(gf^k)'/k!$, and the production matrix $P=[g,f]^{-1}\overline{[g,f]}$ has entries
\[
p_{n,k}=\frac{n!}{k!}\big(c_{n-k}+k\,r_{n-k+1}\big),\qquad c(t)=\frac{g'(\bar f)}{g(\bar f)},\quad r(t)=f'(\bar f),
\]
where $c_i=[t^i]c(t)$, $r_i=[t^i]r(t)$, and coefficients of negative index are zero.
\item[(E3)] More generally, deleting $j$ rows corresponds to the $j$th derivative, and $[g,f]^{-1}U^j[g,f]=P^j$.
\end{itemize}

\begin{example}\label{ex:scales}
The symbol $\big(\frac1{1-x},\frac{x}{1-x}\big)$ gives a different matrix for each scale $c$, since $[x^m]\frac{1}{1-x}\big(\frac{x}{1-x}\big)^k=\binom{m}{k}$. For $c_n=1$ we obtain Pascal's triangle $\big[\binom{m}{k}\big]$, for $c_n=n!$ the matrix $\big[\frac{m!}{k!}\binom{m}{k}\big]$, and for $c_n=(2n)!$ the matrix
\[
\Big(\frac1{1-x},\frac{x}{1-x}\Big)_{(2n)!}=\Big[\frac{(2m)!}{(2k)!}\binom{m}{k}\Big]=\mibb{\begin{array}{cccc}1\\ 2&1\\ 24&24&1\\ 720&1080&90&1\\ \vdots&&&\ddots\end{array}}.
\]
These three matrices are conjugate to each other by diagonal matrices, and they all have the label $\frac{x}{1-x}$ in $\cL$. This is the situation exploited in Sections \ref{sec:EDR} and \ref{sec:EMR}.
\end{example}

\begin{example}\label{ex:expPascal}
The exponential Riordan array $[\ue^t,t]$ is Pascal's triangle $\big[\binom{n}{k}\big]$, because $d_{n,k}=\frac{n!}{k!}[t^n]\ue^tt^k=\binom{n}{k}$. By (E1), $[\ue^t,t]\,\ue^{at}=\ue^{(1+a)t}$, which is the identity $\sum_k\binom{n}{k}a^k=(1+a)^n$. Here $\bar f=t$, $c(t)=1$, and $r(t)=1$, so by (E2) the production matrix has $p_{k,k}=1$ and $p_{k-1,k}=\frac{(k-1)!}{k!}\,k=1$, and all other entries are zero. Hence $P=I+U$, which expresses the recurrence $\binom{n+1}{k}=\binom{n}{k-1}+\binom{n}{k}$.
\end{example}

We next recall the generalized $A$-sequence.

If $(G,F)_c=(d_{n,k})$ and $A(x)=\sum_ja_jx^j$ satisfies $F=x\,A(F)$, i.e., $A(x)=x/\overline F(x)$, then
\begin{equation}\label{genA}
d_{n+1,k+1}=\sum_{j\geq0}\frac{c_{n+1}\,c_{k+j}}{c_{k+1}\,c_n}\,a_j\,d_{n,k+j},\qquad n,k\geq0.
\end{equation}
Conversely, by \eqref{gen.conj} and Rogers' theorem \cite{Rog} (see also \cite{MRSV, HS}), a lower triangular matrix with nonzero diagonal satisfying \eqref{genA} for a sequence with $a_0\not=0$ is a generalized Riordan array with respect to $c$, whose multiplier $F$ is determined by $F=xA(F)$. In particular, two generalized Riordan arrays have the same multiplier if and only if they have the same $A$-sequence.

\begin{example}\label{ex:Aseq}
For Pascal's triangle $\big(\frac1{1-x},\frac{x}{1-x}\big)$ we have $\overline F=\frac{x}{1+x}$, so $A(x)=1+x$. For the Catalan triangle $(C,xC)$, where $C(x)=\frac{1-\sqrt{1-4x}}{2x}$ is the Catalan generating function, we have $\overline F=x(1-x)$, so $A(x)=\frac{1}{1-x}$. For the ballot array $(C,xC^2)$, the identity $xC^2=C-1$ gives $\overline F=\frac{x}{(1+x)^2}$, so $A(x)=(1+x)^2$. These three arrays therefore have pairwise different multipliers, a fact used in Section \ref{sec:inter}.
\end{example}

Finally, we introduce the matrices that appear throughout the paper.

\begin{definition}\label{def:periodic}
Let $\ell\geq2$. A lower triangular matrix $D=(d_{n,j})$ is {\it $\ell$-periodic} if $d_{n,j}=0$ whenever $n\not\equiv j\pmod\ell$. For $\ell=2$ we say {\it checkerboard}. Its {\it compressions} are
\[
D^{(r)}=(d_{\ell m+r,\ \ell k+r})_{m,k\geq0},\qquad 0\leq r\leq\ell-1;
\]
for $\ell=2$ we write $D^\ev=D^{(0)}$ and $D^\od=D^{(1)}$.
\end{definition}

An $\ell$-periodic matrix is determined by its compressions, and any $\ell$ lower triangular matrices $X_0,\ldots,X_{\ell-1}$ are the compressions of a unique $\ell$-periodic matrix, denoted $X_0\oplus\cdots\oplus X_{\ell-1}$. Equivalently, the $\ell$-periodic matrices are those that preserve each summand of $\bK^{\bN_0}=V_0\oplus\cdots\oplus V_{\ell-1}$, where $V_r$ consists of the sequences supported on indices $\equiv r\pmod\ell$.

\begin{lemma}\label{lem:periodic}
If $D$ and $D'$ are $\ell$-periodic, then $DD'$ is $\ell$-periodic and $(DD')^{(r)}=D^{(r)}D'^{(r)}$ for all $r$. If the diagonal of $D$ is nonzero, then $D^{-1}$ is $\ell$-periodic and $(D^{-1})^{(r)}=(D^{(r)})^{-1}$. Moreover, $U^\ell D$ is $\ell$-periodic with compressions $\overline{D^{(r)}}$.
\end{lemma}

\begin{proof}
In $(DD')_{n,j}=\sum_id_{n,i}d'_{i,j}$ a term can be nonzero only if $i\equiv n$ and $i\equiv j\pmod\ell$. Hence $(DD')_{n,j}=0$ for $n\not\equiv j$, and for $n=\ell m+r$, $j=\ell k+r$ only $i=\ell s+r$ contributes, which gives $(DD')^{(r)}=D^{(r)}D'^{(r)}$. For the inverse, $X=\bigoplus_r(D^{(r)})^{-1}$ satisfies $DX=I$. Finally, $(U^\ell D)_{n,j}=d_{n+\ell,j}$, and $(U^\ell D)_{\ell m+r,\ell k+r}=d_{\ell(m+1)+r,\ell k+r}=\overline{D^{(r)}}_{m,k}$.
\end{proof}

For example, if $D$ is a checkerboard matrix with the first rows $1;\ 0,1;\ 2,0,1;\ 0,12,0,1$, then $D^\ev=\left[\begin{smallmatrix}1&0\\2&1\end{smallmatrix}\right]$ and $D^\od=\left[\begin{smallmatrix}1&0\\12&1\end{smallmatrix}\right]$; inverting each block and reassembling gives the first rows $1;\ 0,1;\ -2,0,1;\ 0,-12,0,1$ of $D^{-1}$.

\section{Exponential double Riordan arrays}\label{sec:EDR}

This section is devoted to exponential double Riordan arrays. An exponential double Riordan array is conjugate to an ordinary double Riordan array by a diagonal matrix, so we begin with a brief account of the ordinary double Riordan arrays of Davenport, Shapiro, and Woodson \cite{DSW}. This account introduces the coordinates used throughout the paper, and its fiber product description is the case $c_n=1$ of the exponential results proved below. For ordinary arrays the coefficients may lie in a commutative ring $\bD$, with ``nonzero'' read as ``invertible''.

A double Riordan array $(g;f_1,f_2)$ has $g$ even with $g(0)\not=0$ and $f_1$, $f_2$ odd of order $1$; its columns are
\[
g,\quad gf_1,\quad gf_1f_2,\quad gf_1(f_1f_2),\quad g(f_1f_2)^2,\quad\ldots
\]
These arrays form the double Riordan group $\DR$ under matrix multiplication \cite{DSW}, and the product is usually written as
\begin{equation}\label{DSW-prod}
(g;f_1,f_2)(d;h_1,h_2)=\Big(g\,d(\sqrt{f_1f_2});\ \sqrt{f_1/f_2}\,h_1(\sqrt{f_1f_2}),\ \sqrt{f_2/f_1}\,h_2(\sqrt{f_1f_2})\Big).
\end{equation}

We now introduce coordinates. Since $g$ and $f_1f_2$ are even and $f_1$ is odd, there are unique series $G$, $\Psi$, $\Fh$ with
\begin{equation}\label{coords2}
g(t)=G(t^2),\qquad f_1(t)=t\,\Psi(t^2),\qquad f_1(t)f_2(t)=\Fh(t^2),
\end{equation}
where $G(0)=g(0)\not=0$, $\Psi(0)=f'_1(0)\not=0$, and $\Fh\in\cL$ with $\Fh'(0)=f'_1(0)f'_2(0)$. Conversely, any such triple determines
\begin{equation}\label{coords2-inv}
g(t)=G(t^2),\qquad f_1(t)=t\,\Psi(t^2),\qquad f_2(t)=t\,\frac{\Fh(t^2)/t^2}{\Psi(t^2)}.
\end{equation}
Thus $(g;f_1,f_2)\leftrightarrow(G,\Psi,\Fh)$ is a bijection, and no square root is involved. We call $(G,\Psi,\Fh)$ the {\it coordinates} of the array.

\begin{theorem}\label{thm:DR}
A double Riordan array $D=(g;f_1,f_2)$ is a checkerboard matrix, and its compressions are the Riordan arrays
\begin{equation}\label{DR-compr}
D^\ev=(G,\Fh),\qquad D^\od=(G\Psi,\Fh).
\end{equation}
The map $D\mapsto(D^\ev,D^\od)$ is a group isomorphism
\[
\DR\ \cong\ \cR\times_\cL\cR=\big\{\big((G_0,F),(G_1,F)\big)\big\},
\]
and in coordinates the group law is
\begin{equation}\label{prod-coords}
(G,\Psi,\Fh)\cdot(G',\Psi',\Fh')=\big(G\,G'(\Fh),\ \Psi\,\Psi'(\Fh),\ \Fh'(\Fh)\big),\qquad (G,\Psi,\Fh)^{-1}=\Big(\frac1{G(\Fb)},\ \frac1{\Psi(\Fb)},\ \Fb\Big),
\end{equation}
where $\Fb$ is the compositional inverse of $\Fh$. In particular, $\DR$ is a group over any commutative ring, with no square root needed.
\end{theorem}

\begin{proof}
The $2k$th column is $G(t^2)\Fh(t^2)^k$ and the $(2k+1)$st column is $t\,G(t^2)\Psi(t^2)\Fh(t^2)^k$, which gives the checkerboard property and \eqref{DR-compr}. The rest is the case $c_n=1$ of Theorem \ref{thm:main2} and Corollary \ref{cor:rootfree2} below, whose proofs use only Lemma \ref{lem:periodic}, Proposition \ref{prop:gen}, and \eqref{coords2-inv}.
\end{proof}

When $\sqrt{f_1f_2}$ exists, \eqref{prod-coords} agrees with \eqref{DSW-prod}. Indeed, if $s=\sqrt{f_1f_2}$, then $s^2=\Fh(t^2)$, so $d(s)=G_d(\Fh(t^2))$ and $\sqrt{f_1/f_2}\,h_1(s)=f_1\Psi_h(\Fh(t^2))$, where $(G_d,\Psi_h,\cdot)$ are the coordinates of $(d;h_1,h_2)$. Compressions of double Riordan arrays were studied in \cite{He18}; Theorem \ref{thm:DR} adds that the group law is exactly the fiber product law.

\begin{example}\label{ex:DR}
For $D=\big(\frac1{1-t^2};\frac{t}{1-t^2},t\big)$, the coordinates are $G=\Psi=\frac1{1-x}$ and $\Fh=\frac{x}{1-x}$. The compressions are Pascal's triangle $\big(\frac1{1-x},\frac{x}{1-x}\big)$ and $\big(\frac1{(1-x)^2},\frac{x}{1-x}\big)$, with the same multiplier, and we have
\[
D=\mibb{\begin{array}{cccccccc}
1\\ 0&1\\ 1&0&1\\ 0&2&0&1\\ 1&0&2&0&1\\ 0&3&0&3&0&1\\ 1&0&3&0&3&0&1\\ 0&4&0&6&0&4&0&1\\ \vdots&&&&&&&\ddots
\end{array}},\qquad D^\ev=\mibb{\begin{array}{cccc}1\\1&1\\1&2&1\\1&3&3&1\end{array}},\quad D^\od=\mibb{\begin{array}{cccc}1\\2&1\\3&3&1\\4&6&4&1\end{array}}.
\]
\end{example}

\begin{example}\label{ex:DRsquare}
We compute the square of the array $D$ of Example \ref{ex:DR} by \eqref{prod-coords}, without any square root. Since $G(\Fh)=\Psi(\Fh)=\frac{1-x}{1-2x}$ and $\Fh(\Fh(x))=\frac{x}{1-2x}$, the coordinates of $D^2$ are $\big(\frac{1}{1-2x},\frac{1}{1-2x},\frac{x}{1-2x}\big)$. By \eqref{coords2-inv},
\[
D^2=\Big(\frac{1}{1-2t^2};\ \frac{t}{1-2t^2},\ t\Big)=\mibb{\begin{array}{cccccccc}
1\\ 0&1\\ 2&0&1\\ 0&4&0&1\\ 4&0&4&0&1\\ 0&12&0&6&0&1\\ 8&0&12&0&6&0&1\\ 0&32&0&24&0&8&0&1\\ \vdots&&&&&&&\ddots
\end{array}},
\]
and its compressions $\big(\frac1{1-2x},\frac{x}{1-2x}\big)$ and $\big(\frac1{(1-2x)^2},\frac{x}{1-2x}\big)$ are the squares of the compressions of $D$.
\end{example}

\begin{example}\label{ex:DRinv}
Let $D=\big(\frac1{1-t^2};\frac{t}{1-t^2},-t\big)$. Its coordinates are $G=\Psi=\frac1{1-x}$ and $\Fh=-\frac{x}{1-x}$. Here $f'_1(0)f'_2(0)=-1$, so $\sqrt{f_1f_2}$ does not exist over $\bR$ and \eqref{DSW-prod} cannot be used. By \eqref{prod-coords}, however, $\Fh(\Fh(x))=x$ and $G\cdot G(\Fh)=\frac{1}{1-x}\cdot(1-x)=1$, and likewise for $\Psi$. Hence $D^2=I$, and
\[
D=\mibb{\begin{array}{cccccccc}
1\\ 0&1\\ 1&0&-1\\ 0&2&0&-1\\ 1&0&-2&0&1\\ 0&3&0&-3&0&1\\ 1&0&-3&0&3&0&-1\\ 0&4&0&-6&0&4&0&-1\\ \vdots&&&&&&&\ddots
\end{array}}.
\]
Its compressions $\big(\frac1{1-x},-\frac{x}{1-x}\big)=\big[(-1)^k\binom{m}{k}\big]$ and $\big(\frac1{(1-x)^2},-\frac{x}{1-x}\big)=\big[(-1)^k\binom{m+1}{k+1}\big]$ are the signed Pascal-type involutions.
\end{example}

We now turn to the main subject of this section. An {\it exponential double Riordan array} \cite{He-EDRA} is the matrix
\begin{equation}\label{EDRA.def}
[g;f_1,f_2]=\left(g,\ \frac{gf_{1}}{1!},\ \frac{gf_{1}f_{2}}{2!},\ \frac{gf_{1}^{2}f_{2}}{3!},\ \frac{gf_{1}^{2}f_{2}^{2}}{4!},\ \ldots\right),
\end{equation}
whose columns are given by their exponential generating functions, where $g$ is even with $g(0)\not=0$ and $f_1$, $f_2$ are odd of order $1$. Thus $d_{n,2k}=\frac{n!}{(2k)!}[t^n]g(f_1f_2)^k$ and $d_{n,2k+1}=\frac{n!}{(2k+1)!}[t^n]gf_1(f_1f_2)^k$, and
\begin{equation}\label{EDRA-conj}
[g;f_1,f_2]=\Lambda\,(g;f_1,f_2)\,\Lambda^{-1},\qquad \Lambda={\rm diag}(0!,1!,2!,\ldots).
\end{equation}
In \cite{He-EDRA} it is shown that these arrays form a group $\EDR$ with the product \eqref{DSW-prod} (with square brackets) and the inverse
\begin{equation}\label{inv-old}
[g;f_1,f_2]^{-1}=\left[ \frac{1}{g(\fb)};\ \frac{t\fb}{f_1(\fb)},\ \frac{t\fb}{f_2(\fb)}\right],
\end{equation}
where $\fb$ is the compositional inverse of $\sqrt{f_1f_2}$; moreover $D^{-1}$ acts on even series $B$ and odd series $C$ by
\begin{equation}\label{inv-action}
D^{-1}B=\frac{B(\fb)}{g(\fb)},\qquad D^{-1}C=\frac{t\,C(\fb)}{g(\fb)f_1(\fb)}.
\end{equation}
Over $\bR$ these formulas require $f'_1(0)f'_2(0)>0$. We now remove the square root.

We use the coordinates $(G,\Psi,\Fh)$ of \eqref{coords2}--\eqref{coords2-inv}, now over a field $\bK$ of characteristic zero, and the scales
\[
c^{\ev}_n=(2n)!,\qquad c^{\od}_n=(2n+1)!,\qquad \cR_\ev=\cR_{c^\ev},\quad \cR_\od=\cR_{c^\od},\quad \pi_\ev=\pi_{c^\ev},\quad \pi_\od=\pi_{c^\od}.
\]

\begin{proposition}\label{prop:compr2}
Let $D=[g;f_1,f_2]$ have coordinates $(G,\Psi,\Fh)$. Then $D$ is a checkerboard matrix and
\begin{equation}\label{compr2}
D^\ev=(G,\,\Fh)_{(2n)!},\qquad D^\od=(G\Psi,\,\Fh)_{(2n+1)!}.
\end{equation}
\end{proposition}

\begin{proof}
The $(2k)$th column of $D$ has exponential generating function $G(t^2)\Fh(t^2)^k/(2k)!$, which is even. Hence $d_{2m+1,2k}=0$ and
\[
d_{2m,2k}=\frac{(2m)!}{(2k)!}[t^{2m}]G(t^2)\Fh(t^2)^k=\frac{(2m)!}{(2k)!}[x^m]G(x)\Fh(x)^k,
\]
the $(m,k)$ entry of $(G,\Fh)_{(2n)!}$. The $(2k+1)$st column is $t\,G(t^2)\Psi(t^2)\Fh(t^2)^k/(2k+1)!$, which is odd, so $d_{2m,2k+1}=0$ and $d_{2m+1,2k+1}=\frac{(2m+1)!}{(2k+1)!}[x^m]G\Psi\Fh^k$.
\end{proof}

We can now state the main result of this section.

\begin{theorem}\label{thm:main2}
The map
\[
\Phi:\ \EDR\longrightarrow \cR_\ev\times_{\cL}\cR_\od,\qquad \Phi(D)=(D^\ev,D^\od),
\]
is a group isomorphism. Equivalently, the exponential double Riordan arrays are exactly the matrices $X\oplus Y$ with $X\in\cR_\ev$, $Y\in\cR_\od$, and $\pi_\ev(X)=\pi_\od(Y)$. In particular, $\EDR$ is a group under matrix multiplication over any field of characteristic zero, with no condition on $f'_1(0)f'_2(0)$ other than $f'_1(0)f'_2(0)\not=0$.
\end{theorem}

\begin{proof}
By Proposition \ref{prop:compr2}, $\Phi(D)$ is a pair of generalized Riordan arrays with the same multiplier $\Fh$, so $\Phi$ maps into the fiber product. It is injective because a checkerboard matrix is determined by its compressions. It is surjective. Indeed, given $X=(G_0,\Fh)_{(2n)!}$ and $Y=(G_1,\Fh)_{(2n+1)!}$, the coordinates $(G,\Psi,\Fh)=(G_0,G_1/G_0,\Fh)$ define, by \eqref{coords2-inv}, an array $D$ with $\Phi(D)=(X,Y)$. By Lemma \ref{lem:periodic}, $\Phi$ is compatible with products and inverses of checkerboard matrices with nonzero diagonal. Hence $\EDR=\Phi^{-1}(\cR_\ev\times_\cL\cR_\od)$ is closed under multiplication and inversion, and $\Phi$ is an isomorphism.
\end{proof}

By Proposition \ref{prop:pi}, every $\cR_c$ is isomorphic to $\cR$ compatibly with the label maps, so
\[
\EDR\ \cong\ \cR_\ev\times_\cL\cR_\od\ \cong\ \cR\times_\cL\cR\ \cong\ \DR;
\]
explicitly, the composite isomorphism is $[g;f_1,f_2]\mapsto(g;f_1,f_2)$, i.e., $D\mapsto\Lambda^{-1}D\Lambda$, in agreement with \eqref{EDRA-conj}. The exponential version keeps the factorial normalization, since its compressions are generalized Riordan arrays with respect to $(2n)!$ and $(2n+1)!$.

\begin{corollary}\label{cor:rootfree2}
In coordinates, products and inverses in $\EDR$ are given by \eqref{prod-coords}. Explicitly,
\begin{equation}\label{inv-new-g}
[g;f_1,f_2]^{-1}=\left[\frac{1}{G(\Fb(t^2))};\ \frac{t}{\Psi(\Fb(t^2))},\ \frac{\Fb(t^2)\,\Psi(\Fb(t^2))}{t}\right].
\end{equation}
\end{corollary}

\begin{proof}
By Theorem \ref{thm:main2} and \eqref{gen.prod}, $(DD')^\ev=(G\,G'(\Fh),\Fh'(\Fh))_{(2n)!}$ and $(DD')^\od=(G\Psi\cdot(G'\Psi')(\Fh),\Fh'(\Fh))_{(2n+1)!}$; dividing the first component of the odd compression by that of the even compression gives $\Psi\,\Psi'(\Fh)$. Similarly $(D^{-1})^\ev=(1/G(\Fb),\Fb)_{(2n)!}$ and $(D^{-1})^\od=(1/(G\Psi)(\Fb),\Fb)_{(2n+1)!}$. Formula \eqref{inv-new-g} follows from \eqref{coords2-inv}.
\end{proof}

When $\sqrt{f_1f_2}$ exists, $\fb(t)^2=\Fb(t^2)$ and $f_1(\fb)=\fb\,\Psi(\Fb(t^2))$, so $t\fb/f_1(\fb)=t/\Psi(\Fb(t^2))$, and \eqref{inv-new-g} agrees with \eqref{inv-old}.

\begin{example}\label{ex:Bell2}
Let $D=[\frac{1}{1-t^2};\frac{t}{1-t^2},t]$. Its coordinates are $G=\Psi=\frac{1}{1-x}$ and $\Fh=\frac{x}{1-x}$, so $\Fb=\frac{x}{1+x}$, $G(\Fb)=\Psi(\Fb)=1+x$, and
\[
\left[\frac{1}{1-t^2};\ \frac{t}{1-t^2},\ t\right]^{-1}=\left[\frac{1}{1+t^2};\ \frac{t}{1+t^2},\ t\right].
\]
Here $D^\ev=\big[\frac{(2m)!}{(2k)!}\binom{m}{k}\big]$ and $D^\od=\big[\frac{(2m+1)!}{(2k+1)!}\binom{m+1}{k+1}\big]$, and $D$ and $D^{-1}$ begin
\[
\mibb{\begin{array}{cccccccc}
1\\ 0&1\\ 2&0&1\\ 0&12&0&1\\ 24&0&24&0&1\\ 0&360&0&60&0&1\\ 720&0&1080&0&90&0&1\\ 0&20160&0&5040&0&168&0&1\\ \vdots&&&&&&&\ddots
\end{array}},\]
\[
\mibb{\begin{array}{cccccccc}
1\\ 0&1\\ -2&0&1\\ 0&-12&0&1\\ 24&0&-24&0&1\\ 0&360&0&-60&0&1\\ -720&0&1080&0&-90&0&1\\ 0&-20160&0&5040&0&-168&0&1\\ \vdots&&&&&&&\ddots
\end{array}}.
\]
The inverse is obtained by inverting the two Pascal-type compressions.
\end{example}

The fiber product description reveals the structure of $\EDR$.

\begin{theorem}\label{thm:normal2}
Let $\cN=\{[g;f_1,f_2]\in\EDR:\ f_1f_2=t^2\}=\big\{[g;f_1,t^2/f_1]\big\}$. Then the following statements hold.
\begin{itemize}
\item[(i)] $\cN$ is an abelian normal subgroup of $\EDR$, isomorphic to $(\bK[[x]]^\times)^2$ via $D\mapsto(G,\Psi)$;
\item[(ii)] $\rho:\EDR\to\cL$, $\rho(D)=\Fh$, is a surjective homomorphism with kernel $\cN$;
\item[(iii)] $\cL'=\{[1;t,\Fh(t^2)/t]:\ \Fh\in\cL\}$ is a subgroup mapped isomorphically onto $\cL$ by $\rho$, and
\[
\EDR=\cN\rtimes\cL'\cong(\bK[[x]]^\times)^2\rtimes\cL,
\]
where $\cL$ acts on $\cN$ by substitution, namely $(1,1,\Fh)(G,\Psi,x)(1,1,\Fh)^{-1}=(G(\Fh),\Psi(\Fh),x)$.
\end{itemize}
\end{theorem}

\begin{proof}
By \eqref{prod-coords}, $\rho$ is a surjective homomorphism, and its kernel consists of the elements with $\Fh=x$, i.e., $f_1f_2=t^2$; this is $\cA_\ev\times\cA_\od$ in \eqref{exact}. For $\Fh=\Fh'=x$, \eqref{prod-coords} reads $(G,\Psi,x)\cdot(G',\Psi',x)=(GG',\Psi\Psi',x)$, which proves (i) and (ii). The elements of $\cL'$ have coordinates $(1,1,\Fh)$ and multiply by $\Fh\star\Fh'$. Every element factors as $(G,\Psi,\Fh)=(G,\Psi,x)\cdot(1,1,\Fh)$, and $\cN\cap\cL'=\{I\}$. The conjugation formula follows from \eqref{prod-coords}.
\end{proof}

\begin{example}\label{ex:split2}
The array $D=[\frac1{1-t^2};\frac{t}{1-t^2},t]$ of Example \ref{ex:Bell2} has coordinates $(G,\Psi,\Fh)=\big(\frac1{1-x},\frac1{1-x},\frac{x}{1-x}\big)$. By Theorem \ref{thm:normal2} it factors as $(G,\Psi,x)\cdot(1,1,\Fh)$, that is,
\[
\left[\frac1{1-t^2};\ \frac{t}{1-t^2},\ t\right]=\left[\frac1{1-t^2};\ \frac{t}{1-t^2},\ t(1-t^2)\right]\cdot\left[1;\ t,\ \frac{t}{1-t^2}\right],
\]
where the first factor lies in $\cN$ (indeed $\frac{t}{1-t^2}\cdot t(1-t^2)=t^2$) and the second lies in $\cL'$. The first factor begins with the rows $1;\ 0,1;\ 2,0,1;\ 0,12,0,1;\ 24,0,12,0,1;\ 0,360,0,40,0,1$.
\end{example}

Let $\cA=\{[g;t,t]\}$ (coordinates $(G,1,x)$) and $\cS=\{[1;f_1,f_2]\}$ (coordinates $(1,\Psi,\Fh)$) be the Appell and associated subgroups.

\begin{proposition}\label{prop:subgroups2}
\begin{itemize}
\item[(i)] $\cS\cong\cR_\od\cong\cR$ via $D\mapsto D^\od=(\Psi,\Fh)_{(2n+1)!}$.
\item[(ii)] The diagonal subgroup $\Delta=\{[g;t,f_2]\}$ ($\Psi=1$) is isomorphic to $\cR$; for $D\in\Delta$ both compressions have the symbol $(G,\Fh)$.
\item[(iii)] $\cA\cong\bK[[x]]^\times$ is normal, being the kernel of the homomorphism $(G,\Psi,\Fh)\mapsto(1,\Psi,\Fh)$, and $\EDR=\cA\rtimes\cS\cong\bK[[x]]^\times\rtimes\cR$.
\item[(iv)] For the Bell subgroups $\cB_1=\{f_1=tg\}$ and $\cB_2=\{f_2=tg\}$, we have $D\in\cB_1$ if and only if $\Psi=G$, and $D\in\cB_2$ if and only if $\Fh=xG\Psi$, i.e., if and only if $D^\od=(G\Psi,xG\Psi)_{(2n+1)!}$ is a Bell array.
\end{itemize}
\end{proposition}

\begin{proof}
(i) With $G=G'=1$, \eqref{prod-coords} is $(1,\Psi,\Fh)\cdot(1,\Psi',\Fh')=(1,\Psi\Psi'(\Fh),\Fh'(\Fh))$, the Riordan product of $(\Psi,\Fh)$ and $(\Psi',\Fh')$. (ii) With $\Psi=\Psi'=1$, it is the Riordan product of $(G,\Fh)$ and $(G',\Fh')$. (iii) In a product, $\Psi$ and $\Fh$ do not depend on $G$, so the map is a homomorphism; it is the identity on $\cS$, and $(G,\Psi,\Fh)=(G,1,x)\cdot(1,\Psi,\Fh)$. (iv) $f_1=tg$ means $t\Psi(t^2)=tG(t^2)$, and by \eqref{coords2-inv}, $f_2=tg$ means $\Fh(t^2)=t^2G(t^2)\Psi(t^2)$.
\end{proof}

For example, the array of Example \ref{ex:Bell2} lies in $\cB_1$, and the Pascal--Stanley array $[\frac{1}{1-t^2};t,\frac{t}{1-t^2}]$ lies in $\Delta$, with compressions $\big[\frac{(2m)!}{(2k)!}\binom{m}{k}\big]$ and $\big[\frac{(2m+1)!}{(2k+1)!}\binom{m}{k}\big]$.

We next determine the involutions of $\EDR$.

\begin{theorem}\label{thm:involution}
An exponential double Riordan array $D$ with coordinates $(G,\Psi,\Fh)$ satisfies $D^2=I$ if and only if
\begin{equation}\label{inv-cond}
\Fh(\Fh(x))=x,\qquad G(x)\,G(\Fh(x))=1,\qquad \Psi(x)\,\Psi(\Fh(x))=1.
\end{equation}
If $D^2=I$ and $D\notin\cN$, then $f'_1(0)f'_2(0)=\Fh'(0)=-1$.
\end{theorem}

\begin{proof}
By \eqref{prod-coords}, $D^2$ has coordinates $(G\,G(\Fh),\Psi\,\Psi(\Fh),\Fh(\Fh))$, and $I$ has coordinates $(1,1,x)$. If $\Fh(\Fh(x))=x$, then $\Fh'(0)^2=1$. If $\Fh'(0)=1$ and $\Fh\not=x$, write $\Fh=x+ax^s+O(x^{s+1})$ with $s\geq2$, $a\not=0$; then $\Fh(\Fh(x))=x+2ax^s+O(x^{s+1})\not=x$. Hence $\Fh'(0)=-1$ unless $\Fh=x$.
\end{proof}

Thus every involution outside $\cN$ has $f'_1(0)f'_2(0)<0$, so over $\bR$ it lies outside the square-root framework of \eqref{DSW-prod}--\eqref{inv-old}.

\begin{example}\label{ex:involution}
Take $\Fh=-x$ and $G=\Psi=\ue^{x}$, which satisfy \eqref{inv-cond}. By \eqref{coords2-inv},
\[
J=\left[\ue^{t^2};\ t\ue^{t^2},\ -t\ue^{-t^2}\right]=\mibb{\begin{array}{ccccccc}
1\\ 0&1\\ 2&0&-1\\ 0&12&0&-1\\ 12&0&-12&0&1\\ 0&240&0&-40&0&1\\ 120&0&-180&0&30&0&-1\\ \vdots&&&&&&\ddots
\end{array}},\qquad J^2=I,
\]
with $d_{2m,2k}=(-1)^k\frac{(2m)!}{(2k)!\,(m-k)!}$ and $d_{2m+1,2k+1}=(-1)^k\frac{(2m+1)!}{(2k+1)!}\frac{2^{m-k}}{(m-k)!}$. Its compressions $J^\ev=(\ue^x,-x)_{(2n)!}$ and $J^\od=(\ue^{2x},-x)_{(2n+1)!}$ are involutions in $\cR_\ev$ and $\cR_\od$. More generally, any two involutions $(G_0,\Fh)_{(2n)!}$ and $(G_1,\Fh)_{(2n+1)!}$ with the same $\Fh$ give the involution with coordinates $(G_0,G_1/G_0,\Fh)$.
\end{example}

\begin{example}\label{ex:EDRinv}
The coordinates of Example \ref{ex:DRinv} also define an involution of $\EDR$, namely
\[
\left[\frac1{1-t^2};\ \frac{t}{1-t^2},\ -t\right]=\mibb{\begin{array}{cccccccc}
1\\ 0&1\\ 2&0&-1\\ 0&12&0&-1\\ 24&0&-24&0&1\\ 0&360&0&-60&0&1\\ 720&0&-1080&0&90&0&-1\\ 0&20160&0&-5040&0&168&0&-1\\ \vdots&&&&&&&\ddots
\end{array}}.
\]
Its entries agree, up to the signs $(-1)^{\lfloor k/2\rfloor}$ in column $k$, with those of the array of Example \ref{ex:Bell2}, whose coordinates differ only in the sign of $\Fh$.
\end{example}

We now turn to production matrices.

Let $P=D^{-1}UD$ be the production matrix of $D\in\EDR$. Since $U$ changes parity, $p_{n,j}=0$ whenever $n+j$ is even \cite{He-EDRA}. We describe its two nonzero blocks
\[
P^{\rm eo}=(p_{2m,2k+1})_{m,k\geq0},\qquad P^{\rm oe}=(p_{2m+1,2k})_{m,k\geq0}.
\]

\begin{theorem}\label{thm:P2}
Let $E=D^\ev$ and $O=D^\od$. Then
\begin{equation}\label{blocks2}
P^{\rm eo}=E^{-1}O,\qquad P^{\rm oe}=O^{-1}\overline{E},
\end{equation}
and
\begin{equation}\label{P2}
P^2=P_E\oplus P_O,
\end{equation}
where $P_E=E^{-1}\overline E$ and $P_O=O^{-1}\overline O$ are the production matrices of the compressions. Equivalently, $P^{\rm eo}P^{\rm oe}=P_E$ and $P^{\rm oe}P^{\rm eo}=P_O$.
\end{theorem}

\begin{proof}
Since $(UD)_{n,k}=d_{n+1,k}$, we have $(UD)_{2m,2k+1}=O_{m,k}$, $(UD)_{2m+1,2k}=d_{2m+2,2k}=\overline E_{m,k}$, and all other entries vanish. As $D^{-1}=E^{-1}\oplus O^{-1}$, we get $p_{2m,2k+1}=(E^{-1}O)_{m,k}$ and $p_{2m+1,2k}=(O^{-1}\overline E)_{m,k}$. Further, $P^2=D^{-1}U^2D$, and by Lemma \ref{lem:periodic}, $U^2D$ is a checkerboard matrix with compressions $\overline E$ and $\overline O$, which gives \eqref{P2}.
\end{proof}

We call $P^2$ the {\it second-order production matrix}. Deleting two rows corresponds to $d^2/dt^2$ (fact (E3)), and with \eqref{inv-action} we obtain an explicit formula.

\begin{theorem}\label{thm:P2explicit}
Let $D=[g;f_1,f_2]$ and $F=f_1f_2$, and define the even power series
\begin{align*}
&\alpha(t)=\frac{g''(\fb)}{g(\fb)},\qquad \beta(t)=\frac{(g'F')(\fb)}{g(\fb)},\qquad
\gamma(t)=F''(\fb),\qquad \delta(t)=F'(\fb)^2,\\
&\hat\alpha(t)=\frac{(gf_1)''(\fb)}{(gf_1)(\fb)},\qquad \hat\beta(t)=\frac{((gf_1)'F')(\fb)}{(gf_1)(\fb)}.
\end{align*}
Then the entries $q_{n,j}$ of $P^2$ are, with coefficients of negative index equal to zero,
\begin{align}
q_{n,2k}&=\frac{n!}{(2k)!}\Big(\alpha_{n-2k}+k\,(2\beta+\gamma)_{n-2k+2}+k(k-1)\,\delta_{n-2k+4}\Big),\label{q-even}\\
q_{n,2k+1}&=\frac{n!}{(2k+1)!}\Big(\hat\alpha_{n-2k-1}+k\,(2\hat\beta+\gamma)_{n-2k+1}+k(k-1)\,\delta_{n-2k+3}\Big),\label{q-odd}
\end{align}
where $\alpha_i=[t^i]\alpha(t)$, and similarly for the other series. Along each line $n-2k=$ const, the scaled entries $\frac{(2k)!}{n!}q_{n,2k}$ and $\frac{(2k+1)!}{n!}q_{n,2k+1}$ are quadratic polynomials in $k$.
\end{theorem}

\begin{proof}
The $(2k)$th column of $U^2D$ is $(gF^k)''/(2k)!$, which is even, and
\[
(gF^k)''=g''F^k+2k\,g'F'F^{k-1}+k\,gF''F^{k-1}+k(k-1)\,gF'^2F^{k-2}.
\]
By \eqref{inv-action} and $F(\fb)=t^2$, the $(2k)$th column of $P^2$ is
\[
\frac{1}{(2k)!}\Big(\alpha(t)t^{2k}+k\big(2\beta(t)+\gamma(t)\big)t^{2k-2}+k(k-1)\delta(t)t^{2k-4}\Big),
\]
and $n![t^n]$ gives \eqref{q-even}. The $(2k+1)$st column of $U^2D$ is $(gf_1F^k)''/(2k+1)!$, which is odd; the same expansion with $g$ replaced by $gf_1$ and the odd case of \eqref{inv-action} give \eqref{q-odd}. All six series are even, being quotients of products of series of matching parity composed with the odd series $\fb$.
\end{proof}

All six series are functions of $\fb^{\,2}=\Fb(t^2)$, so they can be computed without square roots; hence Theorem \ref{thm:P2explicit} holds over any field of characteristic zero.

\begin{corollary}\label{cor:twostep}
The entries of $D$ satisfy the two-step recurrences
\[
d_{n+2,2k}=\sum_{i}d_{n,2i}\,q_{2i,2k},\qquad d_{n+2,2k+1}=\sum_{i}d_{n,2i+1}\,q_{2i+1,2k+1},
\]
which involve only one compression each; $(q_{2m,2k})$ and $(q_{2m+1,2k+1})$ are the production matrices of $D^\ev$ and $D^\od$.
\end{corollary}

\begin{example}\label{ex:PS-P2}
For the Pascal--Stanley array $D=[\frac{1}{1-t^2};t,\frac{t}{1-t^2}]$ we have $F=\frac{t^2}{1-t^2}$, $1-\fb^{\,2}=\frac{1}{1+t^2}$, and
\begin{align*}
&\alpha=2(1+t^2)(1+4t^2),\quad \beta=4t^2(1+t^2)^2,\quad \gamma=2(1+t^2)^2(1+4t^2),\quad \delta=4t^2(1+t^2)^3,\\
&\hat\alpha=2(1+t^2)(3+4t^2),\quad \hat\beta=2(1+t^2)^2(1+2t^2).
\end{align*}
The second-order production matrix is a band matrix, which begins
\[
P^2=\mibb{\begin{array}{cccccccc}
2&0&1\\
0&6&0&1\\
20&0&22&0&1\\
0&84&0&34&0&1\\
192&0&528&0&66&0&1\\
0&960&0&1040&0&86&0&1\\
0&0&8640&0&3060&0&134&0\\
0&0&0&20160&0&4788&0&162\\
\vdots&&&&&&&\ddots
\end{array}}.
\]
For instance, $q_{2,0}=2!\,\alpha_2=20$, $q_{3,1}=3!\,\hat\alpha_2=84$, and $d_{6,2}=\sum_i d_{4,2i}q_{2i,2}=24\cdot 1+24\cdot 22+1\cdot 528=1080$.
\end{example}

Finally, since the multiplier of a generalized Riordan array is determined by its $A$-sequence, the fiber product condition can be expressed by a single $A$-sequence.

\begin{theorem}\label{thm:sharedA2}
A lower triangular matrix $D=(d_{n,k})$ with nonzero diagonal is an exponential double Riordan array if and only if it is a checkerboard matrix and there is a sequence $(a_j)_{j\geq0}$ with $a_0\not=0$ that is an $A$-sequence of both compressions, i.e., for all $m,k\geq 0$,
\begin{align*}
d_{2m+2,2k+2}&=\sum_{j\geq 0}\frac{(2m+2)!\,(2k+2j)!}{(2k+2)!\,(2m)!}\,a_j\,d_{2m,2k+2j},\\
d_{2m+3,2k+3}&=\sum_{j\geq 0}\frac{(2m+3)!\,(2k+2j+1)!}{(2k+3)!\,(2m+1)!}\,a_j\,d_{2m+1,2k+2j+1}.
\end{align*}
In this case $\Fh=xA(\Fh)$, $G(x)=\sum_m d_{2m,0}\,x^m/(2m)!$, and $G(x)\Psi(x)=\sum_m d_{2m+1,1}\,x^m/(2m+1)!$.
\end{theorem}

\begin{proof}
If $D\in\EDR$, apply \eqref{genA} to $D^\ev=(G,\Fh)_{(2n)!}$ and $D^\od=(G\Psi,\Fh)_{(2n+1)!}$ with the common $A(x)=x/\Fb(x)$. Conversely, by the converse of \eqref{genA}, $D^\ev\in\cR_\ev$ and $D^\od\in\cR_\od$ with the same multiplier, so $(D^\ev,D^\od)$ lies in the fiber product and $D\in\EDR$ by Theorem \ref{thm:main2}.
\end{proof}

\begin{example}\label{ex:sharedA2}
For the Pascal--Stanley array $[\frac1{1-t^2};t,\frac{t}{1-t^2}]$ we have $\Fh=\frac{x}{1-x}$, so $A(x)=1+x$, i.e., $a_0=a_1=1$. The first relation of Theorem \ref{thm:sharedA2} with $m=2$, $k=0$ gives
\[
d_{6,2}=\frac{6!\,0!}{2!\,4!}\,d_{4,0}+\frac{6!\,2!}{2!\,4!}\,d_{4,2}=15\cdot24+30\cdot24=1080,
\]
and the second relation with $m=2$, $k=0$ gives $d_{7,3}=\frac{7!\,1!}{3!\,5!}\,d_{5,1}+\frac{7!\,3!}{3!\,5!}\,d_{5,3}=7\cdot120+42\cdot40=2520$, in agreement with the entries of the array.
\end{example}

This explains the characterization of \cite{He-EDRA} by one $A$-sequence and two $Z$-sequences. With $x=t^2$, the $A$-sequence there has generating function $t^2/\fb^{\,2}=A(t^2)$, and the $Z_1$- and $Z_2$-sequences are the $Z$-sequences of the ordinary compressions $(G,\Fh)$ and $(G\Psi,\Fh)$. ``One $A$-sequence and two $Z$-sequences'' is the sequence form of ``two Riordan arrays over the same point of $\cL$''. Likewise, the $A_1$-sequence of \cite{He-EDRA} has generating function $f_1(\fb)/\fb=\Psi(\Fb(t^2))$.

\section{Exponential multiple Riordan arrays}\label{sec:EMR}

Throughout this section $\ell\geq2$, and we write $\fr{0}=1$, $\fr{r}=f_1\cdots f_r$, and $F=\fr{\ell}=f_1\cdots f_\ell$. Our main object is the exponential multiple Riordan array. Since it is conjugate to the corresponding ordinary multiple Riordan array by $\Lambda={\rm diag}(n!)$, we treat the ordinary arrays alongside; all results that do not involve derivatives hold for both, with the same proofs.

Let $g\in\bK[[t^\ell]]$ with $g(0)\not=0$, and let $f_1,\ldots,f_\ell\in t\,\bK[[t^\ell]]$ have order $1$. The ({\it ordinary}) {\it multiple Riordan array} $(g;f_1,\ldots,f_\ell)$ \cite{He-MR} and the {\it exponential multiple Riordan array} $[g;f_1,\ldots,f_\ell]$ \cite{He-EMRA} have the columns
\begin{equation}\label{MR.def}
\mbox{column }\ell k+r:\qquad g\,\fr{r}\,F^k\quad\mbox{and}\quad\frac{g\,\fr{r}\,F^k}{(\ell k+r)!},\qquad k\geq0,\ 0\leq r\leq\ell-1,
\end{equation}
given by ordinary and exponential generating functions, respectively; that is, the multipliers $f_1,\ldots,f_\ell$ are used cyclically. For $\ell=2$ these are the arrays of Section \ref{sec:EDR}. As in \eqref{EDRA-conj},
\begin{equation}\label{EMR-conj}
[g;f_1,\ldots,f_\ell]=\Lambda\,(g;f_1,\ldots,f_\ell)\,\Lambda^{-1},\qquad \Lambda={\rm diag}(n!).
\end{equation}
We denote the sets of these arrays by $\MR_\ell$ and $\EMR_\ell$. In \cite{He-MR, He-EMRA} they are shown to be groups with the product
\begin{equation}\label{prod-oldL}
[g;f_1,\ldots,f_\ell][d;h_1,\ldots,h_\ell]=\left[g\,d(\varphi);\ \frac{f_1}{\varphi}h_1(\varphi),\ldots,\frac{f_\ell}{\varphi}h_\ell(\varphi)\right],\qquad \varphi=F^{1/\ell},
\end{equation}
and similarly for $\MR_\ell$. With $\pb$ the compositional inverse of $\varphi$, the inverse of $D=[g;f_1,\ldots,f_\ell]$ acts on a series $B$ of residue $s$ (i.e., $B\in t^s\bK[[t^\ell]]$, $0\leq s\leq\ell-1$) by \cite[Lemma 2.7]{He-EMRA}
\begin{equation}\label{inv-actionL}
D^{-1}B=\frac{t^s\,B(\pb)}{g(\pb)\,\fr{s}(\pb)}.
\end{equation}
These formulas require an $\ell$th root of $f'_1(0)\cdots f'_\ell(0)$.

As in Section \ref{sec:EDR}, we use coordinates. There are unique $G,\Psi_1,\ldots,\Psi_{\ell-1},\Fh\in\bK[[x]]$ with
\begin{equation}\label{coordsL}
g(t)=G(t^\ell),\qquad f_r(t)=t\,\Psi_r(t^\ell)\ (1\leq r\leq\ell-1),\qquad f_1(t)\cdots f_\ell(t)=\Fh(t^\ell),
\end{equation}
where $G(0)\not=0$, $\Psi_r(0)=f'_r(0)\not=0$, and $\Fh\in\cL$ with $\Fh'(0)=f'_1(0)\cdots f'_\ell(0)$. Conversely, with $\Psr{r}=\Psi_1\cdots\Psi_r$ ($\Psr{0}=1$), any such data define
\begin{equation}\label{coordsL-inv}
g(t)=G(t^\ell),\qquad f_r(t)=t\,\Psi_r(t^\ell)\ (1\leq r\leq\ell-1),\qquad f_\ell(t)=t\,\frac{\Fh(t^\ell)/t^\ell}{\Psr{\ell-1}(t^\ell)}.
\end{equation}
We call $(G,\Psi_1,\ldots,\Psi_{\ell-1},\Fh)$ the coordinates of the array; no roots are involved.

For $0\leq r\leq\ell-1$ let
\[
c^{(r)}_n=(\ell n+r)!,\qquad \cR^{(r)}=\cR_{c^{(r)}},\qquad \pi^{(r)}=\pi_{c^{(r)}}.
\]

\begin{proposition}\label{prop:comprL}
Let $D=(g;f_1,\ldots,f_\ell)$ or $D=[g;f_1,\ldots,f_\ell]$ have coordinates $(G,\Psi_1,\ldots,\Psi_{\ell-1},\Fh)$. Then $D$ is $\ell$-periodic, and its compressions are
\begin{equation}\label{comprL}
(g;f_1,\ldots,f_\ell)^{(r)}=\big(G\,\Psr{r},\ \Fh\big),\qquad [g;f_1,\ldots,f_\ell]^{(r)}=\big(G\,\Psr{r},\ \Fh\big)_{(\ell n+r)!},\qquad 0\leq r\leq\ell-1.
\end{equation}
\end{proposition}

\begin{proof}
By \eqref{coordsL}, $g\fr{r}F^k=t^r\,G(t^\ell)\Psr{r}(t^\ell)\Fh(t^\ell)^k$, which contains only powers $t^n$ with $n\equiv r\pmod\ell$. Hence $D$ is $\ell$-periodic, and $[t^{\ell m+r}]\,g\fr{r}F^k=[x^m]\,G\Psr{r}\Fh^k$. In the exponential case the entry is multiplied by $(\ell m+r)!/(\ell k+r)!=c^{(r)}_m/c^{(r)}_k$.
\end{proof}

\begin{theorem}\label{thm:mainL}
The maps $D\mapsto(D^{(0)},\ldots,D^{(\ell-1)})$ are group isomorphisms
\[
\MR_\ell\ \cong\ \cR\times_\cL\cR\times_\cL\cdots\times_\cL\cR,\qquad
\EMR_\ell\ \cong\ \cR^{(0)}\times_\cL\cR^{(1)}\times_\cL\cdots\times_\cL\cR^{(\ell-1)}
\]
($\ell$ factors). Equivalently, the (exponential) multiple Riordan arrays are exactly the matrices $X_0\oplus\cdots\oplus X_{\ell-1}$ with $X_r$ in $\cR$ (respectively $\cR^{(r)}$) and $X_0,\ldots,X_{\ell-1}$ all having the same multiplier. In particular, $\MR_\ell$ is a group over any commutative ring and $\EMR_\ell$ over any field of characteristic zero, with no condition on $f'_1(0)\cdots f'_\ell(0)$ other than being nonzero (invertible), and $\MR_\ell\cong\EMR_\ell$.
\end{theorem}

\begin{proof}
By Proposition \ref{prop:comprL} the map lands in the fiber product, and it is injective because an $\ell$-periodic matrix is determined by its compressions. Given $X_r=(G_r,\Fh)$ (respectively $(G_r,\Fh)_{c^{(r)}}$), the coordinates
\[
G=G_0,\qquad \Psi_r=G_r/G_{r-1}\quad(1\leq r\leq\ell-1),\qquad \Fh
\]
define, by \eqref{coordsL-inv}, an array whose compressions are $X_0,\ldots,X_{\ell-1}$; so the map is surjective. By Lemma \ref{lem:periodic} it is compatible with products and inverses, so the set of arrays is closed under both and the map is an isomorphism. Finally, both fiber products are isomorphic to $\cR\times_\cL\cdots\times_\cL\cR$ by Proposition \ref{prop:pi}.
\end{proof}

For $\ell=2$, Theorem \ref{thm:mainL} contains Theorems \ref{thm:DR} and \ref{thm:main2}.

\begin{corollary}\label{cor:rootfreeL}
In coordinates, the product and inverse in $\MR_\ell$ and $\EMR_\ell$ are
\begin{align}
&(G,\Psi_1,\ldots,\Psi_{\ell-1},\Fh)\cdot(G',\Psi'_1,\ldots,\Psi'_{\ell-1},\Fh')\nonumber\\
&\qquad=\big(G\,G'(\Fh),\ \Psi_1\Psi'_1(\Fh),\ \ldots,\ \Psi_{\ell-1}\Psi'_{\ell-1}(\Fh),\ \Fh'(\Fh)\big),\label{prod-newL}\\
&(G,\Psi_1,\ldots,\Psi_{\ell-1},\Fh)^{-1}=\left(\frac{1}{G(\Fb)},\ \frac{1}{\Psi_1(\Fb)},\ \ldots,\ \frac{1}{\Psi_{\ell-1}(\Fb)},\ \Fb\right).\label{inv-newL}
\end{align}
Explicitly,
\[
[g;f_1,\ldots,f_\ell]^{-1}=\left[\frac{1}{G(\Fb(t^\ell))};\ \frac{t}{\Psi_1(\Fb(t^\ell))},\ \ldots,\ \frac{t}{\Psi_{\ell-1}(\Fb(t^\ell))},\ \frac{\Fb(t^\ell)\,\Psr{\ell-1}(\Fb(t^\ell))}{t^{\ell-1}}\right],
\]
and the same formula holds for $(g;f_1,\ldots,f_\ell)^{-1}$.
\end{corollary}

\begin{proof}
The $r$th compression of $DD'$ is $\big(G\Psr{r}\cdot(G'\Psi'_{[r]})(\Fh),\ \Fh'(\Fh)\big)$ by \eqref{gen.prod}; quotients of consecutive first components give \eqref{prod-newL}. Inverses are handled in the same way.
\end{proof}

When $\varphi=F^{1/\ell}$ exists, $\varphi^\ell=\Fh(t^\ell)$ and $\pb^{\,\ell}=\Fb(t^\ell)$, and Corollary \ref{cor:rootfreeL} agrees with \eqref{prod-oldL} and the inverse formula of \cite{He-EMRA}.

\begin{example}\label{ex:MR3}
(a) The ordinary $3$-tuple array $D=\big(\frac{1}{1-t^3};\frac{t}{1-t^3},\frac{t}{1-t^3},t\big)$ has coordinates $G=\Psi_1=\Psi_2=\frac1{1-x}$ and $\Fh=\frac{x}{(1-x)^2}$. Its three compressions are $\big(\frac{1}{(1-x)^{r+1}},\frac{x}{(1-x)^2}\big)$, $r=0,1,2$, with the same multiplier, and we have
\[
D=\mibb{\begin{array}{cccccccccccc}
1\\ 0&1\\ 0&0&1\\ 1&0&0&1\\ 0&2&0&0&1\\ 0&0&3&0&0&1\\ 1&0&0&3&0&0&1\\ 0&3&0&0&4&0&0&1\\ 0&0&6&0&0&5&0&0&1\\ 1&0&0&6&0&0&5&0&0&1\\ \vdots&&&&&&&&&&\ddots
\end{array}},\]
\[
D^{(0)}=\mibb{\begin{array}{cccc}1\\1&1\\1&3&1\\1&6&5&1\end{array}},\quad
D^{(1)}=\mibb{\begin{array}{cccc}1\\2&1\\3&4&1\\4&10&6&1\end{array}}.
\]
(b) For the exponential Pascal--Stanley-type array $D=[\frac{1}{1-t^\ell};t,\ldots,t,\frac{t}{1-t^\ell}]$ we have $G=\frac1{1-x}$, $\Psi_r=1$, and $\Fh=\frac{x}{1-x}$, so $\Fb=\frac{x}{1+x}$ and
\[
\left[\frac{1}{1-t^\ell};\ t,\ldots,t,\ \frac{t}{1-t^\ell}\right]^{-1}=\left[\frac{1}{1+t^\ell};\ t,\ldots,t,\ \frac{t}{1+t^\ell}\right].
\]
All compressions have the symbol $(\frac{1}{1-x},\frac{x}{1-x})$, i.e., $D^{(r)}=\big[\frac{(\ell m+r)!}{(\ell k+r)!}\binom{m}{k}\big]$. For $\ell=3$,
\[
D=\mibb{\begin{array}{cccccccccc}
1\\ 0&1\\ 0&0&1\\ 6&0&0&1\\ 0&24&0&0&1\\ 0&0&60&0&0&1\\ 720&0&0&240&0&0&1\\ 0&5040&0&0&420&0&0&1\\ 0&0&20160&0&0&672&0&0&1\\ 362880&0&0&181440&0&0&1512&0&0&1\\ \vdots&&&&&&&&&\ddots
\end{array}}.
\]
\end{example}

In the rest of this section $\mathscr{M}_\ell$ denotes either $\MR_\ell$ or $\EMR_\ell$; all statements follow from \eqref{prod-newL} and hold for both, and we write $\mathscr{M}_1=\cR$.

\begin{theorem}\label{thm:normalL}
Let $\cN=\{D\in\mathscr{M}_\ell:\ f_1\cdots f_\ell=t^\ell\}$. Then $\cN\cong(\bK[[x]]^\times)^\ell$ via $D\mapsto(G,\Psi_1,\ldots,\Psi_{\ell-1})$ is an abelian normal subgroup, $\rho(D)=\Fh$ is a surjective homomorphism $\mathscr{M}_\ell\to\cL$ with kernel $\cN$, and the subgroup $\cL'=\{(1;t,\ldots,t,\Fh(t^\ell)/t^{\ell-1})\}$ (coordinates $(1,\ldots,1,\Fh)$) is a complement, so that
\[
\mathscr{M}_\ell=\cN\rtimes\cL'\cong(\bK[[x]]^\times)^\ell\rtimes\cL.
\]
\end{theorem}

\begin{proof}
The proof is the same as that of Theorem \ref{thm:normal2}. For $\Fh=\Fh'=x$, \eqref{prod-newL} is componentwise multiplication; $\rho$ is a homomorphism by \eqref{prod-newL}; and $(G,\Psi_1,\ldots,\Psi_{\ell-1},\Fh)=(G,\Psi_1,\ldots,\Psi_{\ell-1},x)\cdot(1,\ldots,1,\Fh)$ with $\cN\cap\cL'=\{I\}$.
\end{proof}

Let $\cA=\{(G,1,\ldots,1,x)\}=\{[g;t,\ldots,t]\}$ and $\cS=\{(1,\Psi_1,\ldots,\Psi_{\ell-1},\Fh)\}=\{[1;f_1,\ldots,f_\ell]\}$ be the Appell and associated subgroups. Since $\Psi_r$ and $\Fh$ in a product do not depend on $G$, the map $(G,\Psi_1,\ldots,\Fh)\mapsto(1,\Psi_1,\ldots,\Fh)$ is a homomorphism with kernel $\cA$ that is the identity on $\cS$; hence $\mathscr{M}_\ell=\cA\rtimes\cS$ with $\cA\cong\bK[[x]]^\times$.

\begin{theorem}\label{thm:iterated}
The map $D\mapsto(D^{(1)},\ldots,D^{(\ell-1)})$ is an isomorphism from the associated subgroup $\cS$ onto the $(\ell-1)$-fold fiber product of the corresponding (generalized) Riordan groups, and hence $\cS\cong\mathscr{M}_{\ell-1}$. Consequently,
\[
\mathscr{M}_\ell\cong\bK[[x]]^\times\rtimes\mathscr{M}_{\ell-1}\cong\bK[[x]]^\times\rtimes\Big(\bK[[x]]^\times\rtimes\big(\cdots\rtimes(\bK[[x]]^\times\rtimes\cR)\cdots\big)\Big),
\]
with $\ell-1$ factors $\bK[[x]]^\times$. Thus the exponential multiple Riordan group (and likewise the multiple Riordan group) is obtained from the Riordan group by $\ell-1$ successive split extensions by the Appell group.
\end{theorem}

\begin{proof}
For $D\in\cS$, $G=1$, so $D^{(0)}=(1,\Fh)$ is determined by the common multiplier $\Fh$ of $D^{(1)},\ldots,D^{(\ell-1)}$. Hence $D\mapsto(D^{(1)},\ldots,D^{(\ell-1)})$ is an injective homomorphism into the $(\ell-1)$-fold fiber product, and it is surjective since $(G_1,\Fh),\ldots,(G_{\ell-1},\Fh)$ come from the coordinates $(1,G_1,G_2/G_1,\ldots,G_{\ell-1}/G_{\ell-2},\Fh)$. By Theorem \ref{thm:mainL} (for $\ell-1\geq2$) and Proposition \ref{prop:pi}, this fiber product is isomorphic to $\mathscr{M}_{\ell-1}$; for $\ell=2$ it is isomorphic to $\cR$. Induction finishes the proof.
\end{proof}

Equivalently, forgetting the last compression gives a chain of split surjective homomorphisms with Appell kernels,
\[
\cR\ \longleftarrow\ \mathscr{M}_2\ \longleftarrow\ \mathscr{M}_3\ \longleftarrow\ \mathscr{M}_4\ \longleftarrow\ \cdots,
\]
where $p(X_0,\ldots,X_{\ell-1})=(X_0,\ldots,X_{\ell-2})$ and the section $s(X_0,\ldots,X_{\ell-2})=(X_0,\ldots,X_{\ell-2},X_{\ell-2})$ repeats the last block.

\begin{example}\label{ex:iterated}
Let $D=[\frac1{1-t^3};t,t,\frac{t}{1-t^3}]$, with coordinates $\big(\frac1{1-x},1,1,\frac{x}{1-x}\big)$. Its decomposition in $\cA\rtimes\cS$ is
\[
D=\Big[\frac1{1-t^3};\ t,\ t,\ t\Big]\cdot\Big[1;\ t,\ t,\ \frac{t}{1-t^3}\Big],
\]
since $(G,1,1,x)\cdot(1,1,1,\Fh)=(G,1,1,\Fh)$. The second factor lies in $\cS$, and under Theorem \ref{thm:iterated} it corresponds to the pair of compressions $\big((1,\frac{x}{1-x})_{(3n+1)!},(1,\frac{x}{1-x})_{(3n+2)!}\big)$, i.e., to the element of $\EMR_2$ with coordinates $\big(1,1,\frac{x}{1-x}\big)$, which is $[1;t,\frac{t}{1-t^2}]$.
\end{example}

\begin{proposition}\label{prop:BellL}
\begin{itemize}
\item[(i)] The diagonal subgroup $\Delta=\{(g;t,\ldots,t,f_\ell)\}$ ($\Psi_1=\cdots=\Psi_{\ell-1}=1$) is isomorphic to $\cR$; all compressions of $D\in\Delta$ have the same symbol $(G,\Fh)$.
\item[(ii)] For the Bell subgroups $\cB_j=\{f_j=tg\}$, we have $D\in\cB_j$ $(j\leq\ell-1)$ if and only if $\Psi_j=G$, and $D\in\cB_\ell$ if and only if $\Fh=xG\Psr{\ell-1}$, i.e., the last compression $(G\Psr{\ell-1},xG\Psr{\ell-1})$ is a Bell array.
\end{itemize}
\end{proposition}

\begin{proof}
(i) With all $\Psi_r=1$, \eqref{prod-newL} is the Riordan product of $(G,\Fh)$ and $(G',\Fh')$. (ii) $f_j=tg$ means $\Psi_j(t^\ell)=G(t^\ell)$ for $j<\ell$, and, by \eqref{coordsL-inv}, $\Fh(t^\ell)=t^\ell G(t^\ell)\Psr{\ell-1}(t^\ell)$ for $j=\ell$.
\end{proof}

We now turn to elements of finite order.

For $\Fh\in\cL$ let $\Fh^{[i]}$ be its $i$-fold composite ($\Fh^{[0]}=x$).

\begin{theorem}\label{thm:finite}
Let $D\in\mathscr{M}_\ell$ have coordinates $(G,\Psi_1,\ldots,\Psi_{\ell-1},\Fh)$ and $N\geq1$. Then $D^N=I$ if and only if
\begin{equation}\label{finite}
\Fh^{[N]}=x,\qquad \prod_{i=0}^{N-1}G(\Fh^{[i]})=1,\qquad \prod_{i=0}^{N-1}\Psi_r(\Fh^{[i]})=1\quad(1\leq r\leq\ell-1).
\end{equation}
If $D^N=I$, then $\lambda=f'_1(0)\cdots f'_\ell(0)=\Fh'(0)$ satisfies $\lambda^N=1$, and $\lambda=1$ only if $D\in\cN$.
\end{theorem}

\begin{proof}
By induction on \eqref{prod-newL}, $D^N$ has coordinates $\big(\prod_{i<N}G(\Fh^{[i]}),\ \prod_{i<N}\Psi_1(\Fh^{[i]}),\ldots,\Fh^{[N]}\big)$. If $\Fh^{[N]}=x$, then $\Fh'(0)^N=1$; if moreover $\Fh'(0)=1$ and $\Fh=x+ax^s+O(x^{s+1})$ with $s\geq2$, $a\not=0$, then $\Fh^{[N]}=x+Nax^s+O(x^{s+1})\not=x$. Hence $\Fh=x$.
\end{proof}

For $\ell=2$ and $N=2$ this is Theorem \ref{thm:involution}.

\begin{corollary}\label{cor:real}
Over $\bK=\bR$, every element of finite order of $\mathscr{M}_\ell$ is an involution. A nonidentity element of finite order either lies in $\cN$ and has $G,\Psi_1,\ldots,\Psi_{\ell-1}\in\{\pm1\}$, or has $f'_1(0)\cdots f'_\ell(0)=-1$ and satisfies \eqref{finite} with $N=2$.
\end{corollary}

\begin{proof}
If $D\in\cN$ has finite order, then $G^N=1$ in $\bR[[x]]$; writing $G=G(0)(1+h)$ with $h(0)=0$, $(1+h)^N=1$ forces $h=0$, so $G=\pm1$, and similarly for the $\Psi_r$. If $D\notin\cN$, then $\lambda=-1$, and $D^2\in\cN$ has finite order, so its coordinates are constants $\pm1$; but the constant terms of $G\,G(\Fh)$ and $\Psi_r\Psi_r(\Fh)$ are squares, hence positive. So $D^2=I$.
\end{proof}

\begin{example}\label{ex:finite}
(a) For $\ell=3$, $J=[\ue^{t^3};t,t,-t]$ has coordinates $(\ue^x,1,1,-x)$, and $J^2=I$, where
\[
J=\mibb{\begin{array}{cccccccccc}
1\\ 0&1\\ 0&0&1\\ 6&0&0&-1\\ 0&24&0&0&-1\\ 0&0&60&0&0&-1\\ 360&0&0&-120&0&0&1\\ 0&2520&0&0&-210&0&0&1\\ 0&0&10080&0&0&-336&0&0&1\\ 60480&0&0&-30240&0&0&504&0&0&-1\\ \vdots&&&&&&&&&\ddots
\end{array}}.
\]
(b) Over $\bC$, with $\omega=\ue^{2\pi i/3}$, the array $W=[\ue^{t^3};t,t,\omega t]$ has coordinates $(\ue^x,1,1,\omega x)$ and $\prod_{i=0}^2G(\omega^ix)=\ue^{(1+\omega+\omega^2)x}=1$, so $W^3=I\not=W$; for instance $d_{3,3}=\omega$ and $d_{6,3}=120\,\omega$.
\end{example}

\begin{example}\label{ex:MRinv}
Over $\bR$, integer involutions also occur for $\ell=3$. The ordinary array $K=\big(\frac1{1-t^3};t,t,-\frac{t}{1-t^3}\big)$ has coordinates $\big(\frac1{1-x},1,1,-\frac{x}{1-x}\big)$, which satisfy \eqref{finite} with $N=2$. Hence $K^2=I$, and
\[
K=\mibb{\begin{array}{cccccccccc}
1\\ 0&1\\ 0&0&1\\ 1&0&0&-1\\ 0&1&0&0&-1\\ 0&0&1&0&0&-1\\ 1&0&0&-2&0&0&1\\ 0&1&0&0&-2&0&0&1\\ 0&0&1&0&0&-2&0&0&1\\ 1&0&0&-3&0&0&3&0&0&-1\\ \vdots&&&&&&&&&\ddots
\end{array}},
\]
whose three compressions all equal the signed Pascal matrix $\big[(-1)^k\binom{m}{k}\big]$. By \eqref{EMR-conj}, $[\frac1{1-t^3};t,t,-\frac{t}{1-t^3}]=\Lambda K\Lambda^{-1}$ is an involution of $\EMR_3$.
\end{example}

We now turn to production matrices.

Let $D\in\mathscr{M}_\ell$ and $P=D^{-1}UD$. Since $U$ lowers residues by one, $p_{n,j}=0$ unless $n\equiv j-1\pmod\ell$, and $P$ is described by the blocks
\[
P^{(r,r+1)}=(p_{\ell m+r,\ \ell k+r+1})_{m,k\geq0}\ (0\leq r\leq\ell-2),\qquad P^{(\ell-1,0)}=(p_{\ell m+\ell-1,\ \ell k})_{m,k\geq0}.
\]

\begin{theorem}\label{thm:blocks}
We have
\[
P^{(r,r+1)}=\big(D^{(r)}\big)^{-1}D^{(r+1)}\quad(0\leq r\leq\ell-2),\qquad P^{(\ell-1,0)}=\big(D^{(\ell-1)}\big)^{-1}\,\overline{D^{(0)}},
\]
and the $\ell$th power of $P$ is the $\ell$-periodic matrix
\begin{equation}\label{Pl}
P^\ell=P_{D^{(0)}}\oplus P_{D^{(1)}}\oplus\cdots\oplus P_{D^{(\ell-1)}}.
\end{equation}
Moreover, each $P_{D^{(r)}}$ is the cyclic product $P^{(r,r+1)}P^{(r+1,r+2)}\cdots P^{(\ell-1,0)}P^{(0,1)}\cdots P^{(r-1,r)}$.
\end{theorem}

\begin{proof}
$(UD)_{\ell m+r,\ell k+r+1}=D^{(r+1)}_{m,k}$ for $r\leq\ell-2$, $(UD)_{\ell m+\ell-1,\ell k}=d_{\ell(m+1),\ell k}=\overline{D^{(0)}}_{m,k}$, and all other entries of $UD$ vanish; with $D^{-1}=\bigoplus_r(D^{(r)})^{-1}$ this gives the blocks. Next, $P^\ell=D^{-1}U^\ell D$, and Lemma \ref{lem:periodic} gives \eqref{Pl}. The cyclic product telescopes to $(D^{(r)})^{-1}\overline{D^{(0)}}(D^{(0)})^{-1}D^{(r)}=(D^{(r)})^{-1}UD^{(r)}$.
\end{proof}

For $\ell=2$ this is Theorem \ref{thm:P2}. In the exponential case, $U^\ell$ corresponds to $d^\ell/dt^\ell$, and $P^\ell$ can be computed explicitly. Let $B_{j,i}$ denote the partial (exponential) Bell polynomials, so that by Fa\`a di Bruno's formula
\begin{equation}\label{FdB}
\frac{d^j}{dt^j}F^k=\sum_{i=0}^{j}\ff{k}{i}\,F^{k-i}\,B_{j,i}(F',F'',\ldots,F^{(j-i+1)}),
\end{equation}
with $B_{0,0}=1$, $B_{j,0}=0$ for $j\geq1$, and $\ff{k}{i}=k(k-1)\cdots(k-i+1)$.

\begin{theorem}\label{thm:Plexplicit}
Let $D=[g;f_1,\ldots,f_\ell]$, $H_r=g\fr{r}$ $(0\leq r\leq\ell-1)$, and for $0\leq i\leq\ell$ let
\[
\theta_{r,i}(t)=\frac{1}{H_r(\pb)}\sum_{j=i}^{\ell}\binom{\ell}{j}\Big(H_r^{(\ell-j)}\,B_{j,i}(F',F'',\ldots,F^{(j-i+1)})\Big)(\pb).
\]
Then $\theta_{r,i}\in\bK[[t^\ell]]$, and the entries $q_{n,j}$ of $P^\ell$ are
\begin{equation}\label{qformula}
q_{n,\ell k+r}=\frac{n!}{(\ell k+r)!}\sum_{i=0}^{\ell}\ff{k}{i}\,[t^{\,n-r-\ell(k-i)}]\,\theta_{r,i}(t).
\end{equation}
Along each line $n-\ell k=$ const (with $r$ fixed), the scaled entries $\frac{(\ell k+r)!}{n!}q_{n,\ell k+r}$ are polynomials in $k$ of degree at most $\ell$.
\end{theorem}

\begin{proof}
The column $\ell k+r$ of $U^\ell D$ is $(H_rF^k)^{(\ell)}/(\ell k+r)!$, which has residue $r$. By Leibniz's rule and \eqref{FdB},
\[
(H_rF^k)^{(\ell)}=\sum_{i=0}^{\ell}\ff{k}{i}F^{k-i}\sum_{j=i}^{\ell}\binom{\ell}{j}H_r^{(\ell-j)}B_{j,i}(F',\ldots).
\]
By \eqref{inv-actionL} with $s=r$ and $F(\pb)=t^\ell$, the column $\ell k+r$ of $P^\ell$ is $\frac{t^r}{(\ell k+r)!}\sum_i\ff{k}{i}t^{\ell(k-i)}\theta_{r,i}(t)$, which gives \eqref{qformula}. Each monomial of $B_{j,i}(F',\ldots)$ has residue $-j\pmod\ell$ and $H_r^{(\ell-j)}$ has residue $r+j-\ell$, so $\theta_{r,i}$ has residue $0$.
\end{proof}

For $\ell=2$, $\theta_{0,0}=\alpha$, $\theta_{0,1}=2\beta+\gamma$, $\theta_{0,2}=\delta$, and similarly for $r=1$, which is Theorem \ref{thm:P2explicit}. Since $\theta_{r,i}$ is a function of $\pb^{\,\ell}=\Fb(t^\ell)$, it can be computed without roots.

\begin{corollary}\label{cor:lstep}
The entries of $D\in\mathscr{M}_\ell$ satisfy the $\ell$-step recurrences $d_{n+\ell,\ell k+r}=\sum_sd_{n,\ell s+r}\,q_{\ell s+r,\ell k+r}$, each involving only the compression $D^{(r)}$, and $(q_{\ell m+r,\ell k+r})_{m,k}$ is the production matrix of $D^{(r)}$.
\end{corollary}

\begin{example}
For $D=[\frac{1}{1-t^3};t,t,\frac{t}{1-t^3}]$ ($\ell=3$), the third-order production matrix begins
\[
P^3=\mibb{\begin{array}{ccccccccc}
6&0&0&1\\
0&24&0&0&1\\
0&0&60&0&0&1\\
684&0&0&234&0&0&1\\
0&4464&0&0&396&0&0&1\\
0&0&16560&0&0&612&0&0&1\\
194400&0&0&124560&0&0&1272&0&0\\
0&1632960&0&0&282240&0&0&1740&0\\
0&0&7620480&0&0&566496&0&0&2298\\
\vdots&&&&&&&&\ddots
\end{array}},
\]
and its compression on the residue class $0$ is the production matrix of $D^{(0)}=\big[\frac{(3m)!}{(3k)!}\binom{m}{k}\big]$. For example, $d_{9,3}=720\cdot1+240\cdot234+1\cdot124560=181440$, and $q_{0,0}=[t^0]\theta_{0,0}=g'''(0)/g(0)=6$.
\end{example}

We conclude this section with the shared $A$-sequence characterization.

\begin{theorem}\label{thm:sharedAL}
A lower triangular matrix $D=(d_{n,j})$ with nonzero diagonal is an exponential multiple Riordan array if and only if it is $\ell$-periodic and there is a sequence $(a_j)_{j\geq0}$ with $a_0\not=0$ that is an $A$-sequence of every compression, i.e., for $0\leq r\leq\ell-1$ and $m,k\geq0$,
\[
d_{\ell(m+1)+r,\ \ell(k+1)+r}=\sum_{j\geq0}\frac{(\ell(m+1)+r)!\,(\ell(k+j)+r)!}{(\ell(k+1)+r)!\,(\ell m+r)!}\,a_j\,d_{\ell m+r,\ \ell(k+j)+r}.
\]
In this case $\Fh=xA(\Fh)$, and $G\Psr{r}(x)=\sum_md_{\ell m+r,r}\,\frac{r!\,x^m}{(\ell m+r)!}$. For ordinary multiple Riordan arrays the same holds with all factorial quotients replaced by $1$.
\end{theorem}

\begin{proof}
Necessity follows from Proposition \ref{prop:comprL} and \eqref{genA} with the common $A(x)=x/\Fb(x)$. Conversely, by the converse of \eqref{genA}, every compression is a (generalized) Riordan array with the multiplier $\Fh$ determined by $\Fh=xA(\Fh)$; so the compressions lie in the fiber product, and Theorem \ref{thm:mainL} applies.
\end{proof}

\begin{example}\label{ex:sharedAL}
For $D=[\frac1{1-t^3};t,t,\frac{t}{1-t^3}]$ we again have $A(x)=1+x$. The relation of Theorem \ref{thm:sharedAL} with $r=0$, $m=2$, $k=0$ gives
\[
d_{9,3}=\frac{9!\,0!}{3!\,6!}\,d_{6,0}+\frac{9!\,3!}{3!\,6!}\,d_{6,3}=84\cdot720+504\cdot240=181440,
\]
in agreement with Example \ref{ex:MR3}(b).
\end{example}

This explains the characterization of multiple Riordan arrays by one $A$-sequence and $\ell$ $Z$-sequences \cite{He-MR, He-EMRA}. The $A$-sequence has generating function $t^\ell/\pb^{\,\ell}=A(t^\ell)$, and the $Z_r$-sequence is the $Z$-sequence of the ordinary compression $(G\Psr{r},\Fh)$, since $g\fr{r}-\kappa_rt^r=t^\ell g\fr{r}Z_r(\varphi)$ becomes $G\Psr{r}-\kappa_r=x\,G\Psr{r}\,Z_r(\Fh)$, where $\kappa_r=g(0)f'_1(0)\cdots f'_r(0)$. The one-step $A_r$-sequences of \cite{He-EMRA} are $A_r(t)=f_r(\pb)/\pb=\Psi_r(\Fb(t^\ell))$, $1\leq r\leq\ell-1$, and their identity $A_1\cdots A_\ell=A$ expresses $f_1\cdots f_\ell=\Fh(t^\ell)$.

\section{Interleaving Riordan arrays}\label{sec:inter}

Theorems \ref{thm:DR}, \ref{thm:main2}, and \ref{thm:mainL} describe double and multiple Riordan arrays through their compressions. In this section we reverse the point of view. Given Riordan arrays, we ask when they can be interleaved into a double Riordan array, or into an ordinary Riordan array. For lower triangular matrices $X$ and $Y$, the {\it interleaving} $X\oplus Y$ is the checkerboard matrix with compressions $X$ and $Y$ (see Lemma \ref{lem:periodic}).

\begin{theorem}\label{thm:interleave-double}
Let $X=(G_0,F_0)$ and $Y=(G_1,F_1)$ be Riordan arrays. Then $X\oplus Y$ is a double Riordan array if and only if $F_0=F_1$. In that case, writing $F=F_0=F_1$, the double Riordan array is unique and is given by
\begin{equation}\label{double-formula}
g(t)=G_0(t^2),\qquad f_1(t)=t\,\frac{G_1(t^2)}{G_0(t^2)},\qquad f_2(t)=t\,\frac{F(t^2)}{t^2}\cdot\frac{G_0(t^2)}{G_1(t^2)}.
\end{equation}
\end{theorem}

\begin{proof}
If $X\oplus Y=(g;f_1,f_2)$, then by \eqref{DR-compr} both compressions have the multiplier $\Fh$, so $F_0=F_1$. Conversely, let $F_0=F_1=F$ and define $g,f_1,f_2$ by \eqref{double-formula}. Then $g$ is even with $g(0)=G_0(0)\not=0$, $f_1$ is odd with $f'_1(0)=G_1(0)/G_0(0)\not=0$, $f_2$ is odd with $f'_2(0)=F'(0)G_0(0)/G_1(0)\not=0$, and $f_1f_2=F(t^2)$. The coordinates of $(g;f_1,f_2)$ are $(G,\Psi,\Fh)=(G_0,G_1/G_0,F)$, so by \eqref{DR-compr} its compressions are $X$ and $Y$. Uniqueness holds because a checkerboard matrix is determined by its compressions.
\end{proof}

This is the surjectivity of $D\mapsto(D^\ev,D^\od)$ in Theorem \ref{thm:DR}, made explicit. Since the multiplier is determined by the $A$-sequence (see \eqref{genA}), {\it $X$ and $Y$ can be interleaved into a double Riordan array if and only if they have the same $A$-sequence}; their $Z$-sequences may differ.

\begin{remark}\label{rem:types}
Let $X=(G_0,F)$ and $Y=(G_1,F)$. By formula \eqref{double-formula}, the double Riordan array $X\oplus Y$ lies in the diagonal subgroup ($f_1=t$) if and only if $G_1=G_0$. It is a type-$1$ Bell array ($f_1=tg$) if and only if $G_1=G_0^2$. It is a type-$2$ Bell array ($f_2=tg$) if and only if $F=xG_1$, i.e., if and only if $Y=(G_1,xG_1)$ is a Bell array. Finally, it lies in the associated subgroup ($g=1$) if and only if $G_0=1$, i.e., if and only if $X=(1,F)$.
Since the nonzero entries of $X\oplus Y$ are those of $X$ and $Y$, the result is an integer matrix whenever $X$ and $Y$ are, even if $f_1$ and $f_2$ have non-integer coefficients.
\end{remark}

We now ask when the interleaving of two Riordan arrays is an ordinary Riordan array.

\begin{lemma}\label{lem:checkerRiordan}
A Riordan array $(g,f)$ is a checkerboard matrix if and only if $g$ is even and $f$ is odd. In that case, writing $g(t)=G(t^2)$ and $f(t)=t\,\Phi(t^2)$,
\begin{equation}\label{riordan-compr}
(g,f)^\ev=(G,\ x\Phi^2),\qquad (g,f)^\od=(G\Phi,\ x\Phi^2).
\end{equation}
\end{lemma}

\begin{proof}
Column $0$ is $g$, so a checkerboard $(g,f)$ has $g$ even; column $1$ is $gf$, which must be odd, so $f=(gf)/g$ is odd. Conversely, if $g$ is even and $f$ is odd, column $k$, namely $gf^k$, has the parity of $k$. Finally, $gf^{2k}=G(t^2)\big(t^2\Phi(t^2)^2\big)^k$ and $gf^{2k+1}=t\,G(t^2)\Phi(t^2)\big(t^2\Phi(t^2)^2\big)^k$.
\end{proof}

Thus a checkerboard Riordan array is exactly a double Riordan array with $f_1=f_2$.

\begin{theorem}\label{thm:interleave-riordan}
Let $X=(G_0,F_0)$ and $Y=(G_1,F_1)$ be Riordan arrays. Then $X\oplus Y$ is a Riordan array if and only if
\begin{equation}\label{riordan-cond}
F_0=F_1=F\qquad\mbox{and}\qquad F=x\left(\frac{G_1}{G_0}\right)^2.
\end{equation}
In that case
\begin{equation}\label{riordan-formula}
X\oplus Y=(g,f),\qquad g(t)=G_0(t^2),\qquad f(t)=t\,\frac{G_1(t^2)}{G_0(t^2)}.
\end{equation}
\end{theorem}

\begin{proof}
If $X\oplus Y=(g,f)$, then by Lemma \ref{lem:checkerRiordan}, $X=(G,x\Phi^2)$ and $Y=(G\Phi,x\Phi^2)$, so $F_0=F_1=x\Phi^2$ and $G_1/G_0=\Phi$. Conversely, assume \eqref{riordan-cond} and put $\Phi=G_1/G_0$; then \eqref{riordan-formula} defines $(g,f)$ with $g$ even and $f=t\Phi(t^2)$ odd, and by Lemma \ref{lem:checkerRiordan} its compressions are $(G_0,x\Phi^2)=X$ and $(G_0\Phi,x\Phi^2)=Y$. Equivalently, in Theorem \ref{thm:interleave-double} the condition $f_1=f_2$ reads $(G_1/G_0)^2=F/x$.
\end{proof}

\begin{corollary}\label{cor:determined}
Let $X=(G_0,F)$ be given. A Riordan array $Y$ with $X\oplus Y\in\cR$ exists if and only if $F(x)/x$ has a square root in $\bK[[x]]$ (over $\bR$, if and only if $F'(0)>0$), and then $Y=\big(\pm G_0\sqrt{F/x},\ F\big)$; the two signs give $(g,\pm f)$. By contrast, in Theorem \ref{thm:interleave-double} the first function $G_1$ of $Y$ is free.
\end{corollary}

\begin{proposition}\label{prop:checkersubgroup}
The checkerboard Riordan arrays form a subgroup of $\cR$. In the coordinates of $\DR\cong\cR\times_\cL\cR$ it is $\{(G,\Psi,\Fh):\ \Fh=x\Psi^2\}$, and $(g,f)\mapsto\big((G,x\Phi^2),(G\Phi,x\Phi^2)\big)$ is an isomorphism onto it.
\end{proposition}

\begin{proof}
By \eqref{prod-coords}, if $\Fh=x\Psi^2$ and $\Fh'=x\Psi'^2$, then the product has third coordinate $\Fh'(\Fh)=\Fh\,\Psi'(\Fh)^2=x\big(\Psi\,\Psi'(\Fh)\big)^2$, which is $x$ times the square of its second coordinate; inverses are similar. By Lemma \ref{lem:checkerRiordan} these elements are exactly the checkerboard Riordan arrays.
\end{proof}

Both results extend to exponential arrays and to more than two arrays.

An exponential Riordan array $[g,f]$ with $g$ even and $f$ odd is a checkerboard matrix, and the computation in Lemma \ref{lem:checkerRiordan} gives $[g,f]^\ev=(G,x\Phi^2)_{(2n)!}$ and $[g,f]^\od=(G\Phi,x\Phi^2)_{(2n+1)!}$. Together with Proposition \ref{prop:compr2}, this shows that Theorems \ref{thm:interleave-double} and \ref{thm:interleave-riordan} hold verbatim for $X=(G_0,F)_{(2n)!}$ and $Y=(G_1,F)_{(2n+1)!}$, with the same formulas, now producing $[g;f_1,f_2]$ and $[g,f]$. What is interleaved is the pair of symbols $(G_0,F)$ and $(G_1,F)$; the scales are dictated by the position. In particular, an exponential Riordan array with respect to $n!$ must first be converted to its symbol before it can be placed in one of these positions.

For $\ell\geq2$ arrays $X_r=(G_r,F_r)$, $0\leq r\leq\ell-1$, placed on the rows and columns congruent to $r$ modulo $\ell$, Theorem \ref{thm:mainL} and the computation of Lemma \ref{lem:checkerRiordan} give the following two statements.
\begin{itemize}
\item the result is a multiple Riordan array if and only if $F_0=\cdots=F_{\ell-1}=F$, and then
\[
g=G_0(t^\ell),\qquad f_r=t\,\frac{G_r(t^\ell)}{G_{r-1}(t^\ell)}\ (1\leq r\leq\ell-1),\qquad f_\ell=t\,\frac{F(t^\ell)}{t^\ell}\cdot\frac{G_0(t^\ell)}{G_{\ell-1}(t^\ell)};
\]
\item the result is a Riordan array $(g,f)$ if and only if, in addition, $G_r=G_0\Phi^r$ for $0\leq r\leq\ell-1$ and $F=x\Phi^\ell$, where $\Phi=G_1/G_0$; then $g=G_0(t^\ell)$ and $f=t\,\Phi(t^\ell)$, since the column $\ell k+r$ of $(g,f)$ is $t^r\,G(t^\ell)\Phi(t^\ell)^r\big(t^\ell\Phi(t^\ell)^\ell\big)^k$.
\end{itemize}
The exponential versions use the scales $((\ell n+r)!)$ for the block $r$.

We conclude this section with examples.

Let $C(x)=\frac{1-\sqrt{1-4x}}{2x}=1+x+2x^2+5x^3+14x^4+\cdots$ be the Catalan generating function.

\begin{example}\label{ex:pascal}
Pascal's triangle $X=\big(\frac1{1-x},\frac{x}{1-x}\big)$ and $Y=\big(\frac1{(1-x)^2},\frac{x}{1-x}\big)$ have the same multiplier, and \eqref{double-formula} gives the double Riordan array $\big(\frac1{1-t^2};\frac{t}{1-t^2},t\big)$ of Example \ref{ex:DR}. Since $G_1=G_0^2$, it is a type-$1$ Bell array. It is not a Riordan array, since $x(G_1/G_0)^2=\frac{x}{(1-x)^2}\not=\frac{x}{1-x}$.
\end{example}

\begin{example}\label{ex:catalan-double}
$X=(C,xC)$ (the Catalan triangle) and $Y=(C^2,xC)$ have the same multiplier $xC$. By \eqref{double-formula}, $g=C(t^2)$, $f_1=t\,C(t^2)=t+t^3+2t^5+5t^7+\cdots$, $f_2=t$, and
\[
X\oplus Y=\mibb{\begin{array}{cccccccccc}
1\\ 0&1\\ 1&0&1\\ 0&2&0&1\\ 2&0&2&0&1\\ 0&5&0&3&0&1\\ 5&0&5&0&3&0&1\\ 0&14&0&9&0&4&0&1\\ 14&0&14&0&9&0&4&0&1\\ 0&42&0&28&0&14&0&5&0&1\\ \vdots&&&&&&&&&\ddots
\end{array}}.
\]
The compressions are $(C,xC)$, with rows $1;\ 1,1;\ 2,2,1;\ 5,5,3,1;\ldots$, and $(C^2,xC)$, with rows $1;\ 2,1;\ 5,3,1;\ 14,9,4,1;\ldots$. Again $G_1=G_0^2$, so $f_1=tg$; it is not a Riordan array, since $xC^2\not=xC$.
\end{example}

\begin{example}\label{ex:counter}
Pascal's triangle $\big(\frac1{1-x},\frac{x}{1-x}\big)$ and the Catalan triangle $(C,xC)$ are integer Riordan arrays, and their interleaving is a checkerboard matrix. However, their $A$-sequences are $1+x$ and $\frac{1}{1-x}=1+x+x^2+\cdots$, so their multipliers differ, and by Theorem \ref{thm:interleave-double} the interleaving is not a double Riordan array.
\end{example}

\begin{example}\label{ex:mv}
Let $X=\big(\frac1{1-x},\frac{x}{(1-x)^2}\big)$ and $Y=\big(\frac1{(1-x)^2},\frac{x}{(1-x)^2}\big)$, with entries $\binom{m+k}{2k}$ and $\binom{m+k+1}{2k+1}$ (the coefficient triangles of the Morgan--Voyce polynomials $b_m$ and $B_m$). Here $G_1/G_0=\frac1{1-x}$ and $x(G_1/G_0)^2=\frac{x}{(1-x)^2}=F$, so \eqref{riordan-cond} holds, and
\[
X\oplus Y=\Big(\frac1{1-t^2},\ \frac{t}{1-t^2}\Big)=\mibb{\begin{array}{cccccccccc}
1\\ 0&1\\ 1&0&1\\ 0&2&0&1\\ 1&0&3&0&1\\ 0&3&0&4&0&1\\ 1&0&6&0&5&0&1\\ 0&4&0&10&0&6&0&1\\ 1&0&10&0&15&0&7&0&1\\ 0&5&0&20&0&21&0&8&0&1\\ \vdots&&&&&&&&&\ddots
\end{array}},
\]
with $d_{n,k}=\binom{(n+k)/2}{k}$ for $n\equiv k\pmod 2$. Its inverse $\big(\frac1{1+t^2},\frac{t}{1+t^2}\big)$ is the coefficient array of the Chebyshev polynomials $U_n(x/2)$.
\end{example}

\begin{example}\label{ex:ballot}
Let $X=(C,xC^2)$ and $Y=(C^2,xC^2)$. Then $G_1/G_0=C$ and $xC^2=F$, so \eqref{riordan-cond} holds, and
\[
X\oplus Y=\big(C(t^2),\ t\,C(t^2)\big)=\mibb{\begin{array}{cccccccc}
1\\ 0&1\\ 1&0&1\\ 0&2&0&1\\ 2&0&3&0&1\\ 0&5&0&4&0&1\\ 5&0&9&0&5&0&1\\ 0&14&0&14&0&6&0&1\\ \vdots&&&&&&&\ddots
\end{array}},
\]
which is the ballot-number array, i.e., $d_{n,k}$ counts paths of $n$ up and down steps from height $0$ to height $k$ that never go below $0$. Its compressions have rows $1;\ 1,1;\ 2,3,1;\ 5,9,5,1;\ldots$ and $1;\ 2,1;\ 5,4,1;\ 14,14,6,1;\ldots$. Compared with Example \ref{ex:catalan-double}, the odd blocks have the same first function $C^2$; only the multiplier ($xC^2$ versus $xC$) decides whether a Riordan array or only a double Riordan array results.
\end{example}

\begin{example}\label{ex:exp}
The symbols of Example \ref{ex:mv} with the scales $(2n)!$ and $(2n+1)!$ give the exponential Riordan array
\[
\Big[\frac1{1-t^2},\ \frac{t}{1-t^2}\Big]=\mibb{\begin{array}{cccccccc}
1\\ 0&1\\ 2&0&1\\ 0&12&0&1\\ 24&0&36&0&1\\ 0&360&0&80&0&1\\ 720&0&2160&0&150&0&1\\ 0&20160&0&8400&0&252&0&1\\ \vdots&&&&&&&\ddots
\end{array}},
\]
whose compressions have entries $\frac{(2m)!}{(2k)!}\binom{m+k}{2k}$ and $\frac{(2m+1)!}{(2k+1)!}\binom{m+k+1}{2k+1}$.
\end{example}

\begin{example}\label{ex:l3}
(a) $X_r=\big(\frac{1}{(1-x)^{r+1}},\frac{x}{(1-x)^2}\big)$, $r=0,1,2$, share the multiplier and interleave into the multiple Riordan array $\big(\frac{1}{1-t^3};\frac{t}{1-t^3},\frac{t}{1-t^3},t\big)$ of Example \ref{ex:MR3}(a).

(b) $X_r=\big(\frac{1}{(1-x)^{r+1}},\frac{x}{(1-x)^3}\big)$, $r=0,1,2$, satisfy $G_r=G_0\Phi^r$ and $F=x\Phi^3$ with $\Phi=\frac1{1-x}$, so they interleave into the Riordan array $\big(\frac1{1-t^3},\frac{t}{1-t^3}\big)$, whose three compressions have entries $\binom{m+2k+r}{3k+r}$ and begin
\[
\mibb{\begin{array}{cccc}1\\1&1\\1&4&1\\1&10&7&1\end{array}},\qquad
\mibb{\begin{array}{cccc}1\\2&1\\3&5&1\\4&15&8&1\end{array}},\qquad
\mibb{\begin{array}{cccc}1\\3&1\\6&6&1\\10&21&9&1\end{array}}.
\]
\end{example}

In summary, two Riordan arrays $X=(G_0,F_0)$ and $Y=(G_1,F_1)$ can be interleaved into a double Riordan array if and only if $F_0=F_1=F$. In this case the result is
\[
\Big(G_0(t^2);\ t\,\frac{G_1(t^2)}{G_0(t^2)},\ t\,\frac{F(t^2)}{t^2}\cdot\frac{G_0(t^2)}{G_1(t^2)}\Big),
\]
and the first function $G_1$ of $Y$ may be chosen freely. They can be interleaved into an ordinary Riordan array if and only if, in addition, $F=x(G_1/G_0)^2$. In this case the result is $\big(G_0(t^2),\ t\,G_1(t^2)/G_0(t^2)\big)$, and $Y$ is determined by $X$ up to sign. The same statements hold for exponential arrays with the scales $(2n)!$ and $(2n+1)!$, and, as described above, for $\ell$ arrays.

\section{Questions}\label{sec:questions}

\begin{question}
A double or multiple Riordan array with $d_{\ell-1,\ell-1}d_{\ell,0}>0$ is never totally positive, since the minor with rows $\{\ell-1,\ell\}$ and columns $\{0,\ell-1\}$ equals $-d_{\ell-1,\ell-1}d_{\ell,0}$. The natural positivity notion is total positivity of all compressions $(G\Psr{r},\Fh)$. For which coordinates is this the case, and can it be read off from the production matrices of Corollary \ref{cor:lstep}? The total positivity of compressions of ordinary multiple Riordan arrays is studied in \cite{He-MR}, and that of double almost-Riordan arrays in \cite{He-DARA}.
\end{question}

\begin{question}
Classify the conjugacy classes of elements of finite order in $\EMR_\ell$ over $\bC$, and relate them to those in the generalized Riordan groups $\cR^{(r)}$. Which finite groups embed in $\EMR_\ell$? By Theorem \ref{thm:iterated}, this is related to the same question for the Riordan group and its split extensions by Appell groups.
\end{question}

\begin{question}
Riordan pseudo-involutions (arrays $D$ with $(DM)^2=I$ for $M={\rm diag}(1,-1,1,-1,\ldots)$) play a role in the study of the Riordan group. Since $M$ is a checkerboard matrix with $M^\ev=I$ and $M^\od=-I$, we have $(DM)^\ev=D^\ev$ and $(DM)^\od=-D^\od$. What is the corresponding theory for double and exponential double Riordan arrays, and what replaces $M$ for $\ell\geq3$?
\end{question}

\begin{question}
Double almost-Riordan arrays form a group \cite{He-DAR} and have compressions \cite{He-DARA}. Does the fiber product description extend to them, to multiple almost-Riordan arrays, and to the multiple Riordan type arrays of \cite{He-MR}, where a fiber product over a semigroup would replace the fiber product over $\cL$?
\end{question}


\begin{thebibliography}{99}

\bibitem{Barry}
P. Barry, The triple Riordan group, arXiv:2412.05461, 2024.

\bibitem{BDFJT}
D. Branch, D. Davenport, S. Frankson, J. Jones, and G. Thorpe, A and Z sequences for double Riordan arrays, {\it Springer Proceedings in Mathematics and Statistics} 388 (2022), 33--46.

\bibitem{DFSW} D. E. Davenport, S. K. Frankson, L. W. Shapiro, and L. C. Woodson, An invitation to the Riordan group, {\it Enumerative Combinatorics and Applications} 4:3 (2024), Article \#S2S1.

\bibitem{DSW}
D. E. Davenport, L. W. Shapiro, and L. C. Woodson, The double Riordan group, {\it Electronic J. Combin.} 18(2) (2012), P33.

\bibitem{DFR}
E. Deutsch, L. Ferrari, and S. Rinaldi, Production matrices and Riordan arrays, {\it Ann. Combin.} 13 (2009), 65--85.

\bibitem{He15}
T.-X. He, Matrix characterizations of Riordan arrays, {\it Linear Algebra Appl.} 465 (2015), 15--42.

\bibitem{He18}
T.-X. He, Sequence characterizations of double Riordan arrays and their compressions, {\it Linear Algebra Appl.} 549 (2018), 176--202.

\bibitem{He-DAR}
T.-X. He, The double almost-Riordan group, {\it Linear Algebra Appl.} 705 (2025), 50--88.

\bibitem{He-DARA}
T.-X. He, The double almost-Riordan arrays and their sequence characterization, compression, and total positivity, {\it Enumer. Comb. Appl.} 6 (2026), no. 4, Paper No. S2R26, 14 pp.

\bibitem{He-MR}
T.-X. He, The multiple Riordan group and the multiple Riordan type arrays, arXiv:2504.04049, 2025.

\bibitem{He-EDRA}
T.-X. He, Exponential double Riordan arrays, manuscript, 2024.

\bibitem{He-EMRA}
T.-X. He, Exponential multiple Riordan arrays, manuscript, arXiv:2611.15801 [math.CO], 2026.

\bibitem{HS}
T.-X. He and R. Sprugnoli, Sequence characterization of Riordan arrays, {\it Discrete Math.} 309 (2009), 3962--3974.

\bibitem{MRSV}
D. Merlini, D. G. Rogers, R. Sprugnoli, and M. C. Verri, On some alternative characterizations of Riordan arrays, {\it Canadian J. Math.} 49 (1997), 301--320.

\bibitem{Rog}
D. G. Rogers, Pascal triangles, Catalan numbers and renewal arrays, {\it Discrete Math.} 22 (1978), 301--310.

\bibitem{SGWW}
L. W. Shapiro, S. Getu, W. J. Woan, and L. Woodson, The Riordan group, {\it Discrete Appl. Math.} 34 (1991), 229--239.

\end{thebibliography}
\end{document}